\documentclass[12pt]{amsart}

\RequirePackage[dvipsnames,usenames]{color}
\usepackage{hyperref}
\usepackage{amssymb}
\usepackage[latin1]{inputenc}
\usepackage{amsthm}
\usepackage{amsfonts}
\usepackage{mathrsfs}
\usepackage{graphicx}
\usepackage{amsmath}
\usepackage[all,cmtip]{xy}
\numberwithin{equation}{section}

\newtheorem{teo}{Theorem}[section] 
\newtheorem{lem}[teo]{Lemma}

\newtheorem{cor}[teo]{Corollary}
\newtheorem{prop}[teo]{Proposition}

\theoremstyle{definition}
\newtheorem{defi}[teo]{Definition}
\newtheorem{oss}[teo]{Remark}
\newtheorem{ex}[teo]{Example}

\title{On invariance of plurigenera for foliations on surfaces}
\thanks{ Both the authors were partially supported by an EPSRC grant. We would like to thank Y. Gongyo, J. V. Pereira, M. McQuillan, R. Pignatelli and H. Tanaka for several very useful discussions and comments. The work started as an attempt to answer  a question raised by J. V. Pereira during the workshop "Foliation theory in algebraic geometry" supported by the Simons foundation. We would like to thank the referee for carefully reading our manuscript,  for suggesting several improvements and for pointing out  some errors in an earlier version of this paper.  
}

\author{Paolo Cascini} 
\address{Department of Mathematics\\
Imperial College London\\
180 Queen's Gate\\
London SW7 2AZ, UK\\}
\email{p.cascini@imperial.ac.uk}
\author{Enrica Floris}

\email{e.floris@imperial.ac.uk}

\begin{document}
\begin{abstract}
We show that if $(X_t,\mathcal F_t)_{t\in \Delta}$ is a family of foliations  with reduced singularities on a smooth family of surfaces, then invariance of  plurigenera $h^0(X_t,mK_{\mathcal F_t})$ holds for sufficiently large $m$. 
On the other hand, we provide examples on which the result fails, for small values of $m$. 
\end{abstract}
\maketitle

\tableofcontents

\section{Introduction}
The aim of this paper is to study invariance of plurigenera for foliations on algebraic surfaces. 

Both the study of the plurigenera of a manifold and the theory of foliations play a major role in birational geometry. Indeed, on the one hand Siu's Theorem on invariance of plurigenera \cite{siu98,siu02} represents one of the most celebrated results in higher dimensional geometry. The result states that if $\mathcal X\to \Delta$ is a smooth family of projective manifolds over the disk $\Delta$ with fibre $X_t$ at $t\in \Delta$, then the plurigenera $h^0(X_t,\mathcal O_{X_t}(mK_{X_t}))$ does not depend on $t\in \Delta$ (see also \cite{paun05}). 
This generalizes the well-known fact that the genus of a smooth curve is constant under smooth deformations. 
Apart of its own interest, Siu's proof introduces new methods, such as new extension theorems, which had a major impact on some of the recent developments in birational geometry. 

On the other hand, thanks to the work of Miyaoka \cite{miyaoka87}, the theory of foliations plays an important role in the Minimal Model Program in dimension three, as it is needed to  solve some of the crucial cases of the abundance conjecture for threefolds (see also \cite[Chapt. 9]{kollaretal}). It is therefore natural to ask whether the classical results in birational geometry, such as invariance of plurigenera, hold in the more general theory of foliations. 

If $X$ is a smooth surface, a (singular) foliation $\mathcal{F}$ on $X$ corresponds to 
a saturated invertible subsheaf  $T_{\mathcal{F}}\subseteq T_X$.
Its canonical divisor $K_{\mathcal{F}}$ is the divisor associated to the dual of $T_{\mathcal{F}}$.
If $\mathcal{F}$ has reduced singularities 
then the Kodaira dimension $\kappa(\mathcal{F})$
coincides with the Kodaira dimension of its canonical divisor $\kappa(K_{\mathcal{F}})$ 
(see section \ref{s_preliminary} for more details).

Remarkably,  McQuillan \cite{mcq08} was able to reproduce some of the main results of the minimal model program for surfaces to the case of foliations. Similarly, Brunella \cite{brunella01}  showed that
if $(X_t,\mathcal{F}_t)_{t\in\Delta}$ is a family of foliations on surfaces with reduced singularities, then the Kodaira dimension $\kappa(\mathcal{F}_t)$ does not depend on $t$. The following step is therefore to understand to what extent the dimension $h^0(X_t,mK_{\mathcal{F}_t})$  depends on $t\in \Delta$.

\medskip

The goal of this paper is to provide an answer to the question above. More specifically, we prove:

 \medskip

\begin{teo}\label{t_main}
Let $(X_t,\mathcal{F}_t)_{t\in\Delta}$ be a family of foliations with reduced singularities. 
\begin{enumerate}
\item If $\kappa(\mathcal{F}_t)=0$ for any $t\in \Delta$,
then for any $m\in\mathbb{N}$ the dimension $h^0(X_t,\mathcal O_{X_t}(mK_{\mathcal{F}_t}))$ does not depend on $t\in \Delta$;
\item If $\kappa(\mathcal{F}_t)=1$
and $\mathcal{F}_t$ is induced by an elliptic fibration 
then for any sufficiently large positive integer $m$ the dimension $h^0(X_t,\mathcal O_{X_t}(mK_{\mathcal{F}_t}))$
does not depend on $t\in \Delta$;
\item  If $\kappa(\mathcal{F}_t)=1$ and $\mathcal{F}_t$ is not induced by an elliptic fibration for any $t\in \Delta$
then for any positive integer $m$ the dimension $h^0(X_t,\mathcal O_{X_t}(mK_{\mathcal{F}_t}))$ does not depend on $t\in \Delta$; and
\item  If $\kappa(\mathcal F_t)=2$ for any $t\in \Delta$, then 
for any sufficiently large positive integer $m$ the dimension $h^0(X_t,\mathcal O_{X_t}(mK_{\mathcal{F}_t}))$
does not depend on $t\in \Delta$.
\end{enumerate}
\end{teo}

Note that by Brunella's result above, for each family of foliations $(X_t,\mathcal{F}_t)_{t\in\Delta}$ with non-negative Kodaira dimension, either $\kappa(\mathcal F_t)=0$ for any $t\in \Delta$, or $\kappa(\mathcal F_t)=1$ for any $t\in \Delta$ or $\kappa(\mathcal F_t)=2$ for any $t\in \Delta$.
Furthermore we show that, in cases (2) and (4), the invariance of $h^0(X_t,mK_{\mathcal F_t})$ fails for small values of $m$. Indeed, we provide examples of families of foliations  such that
$h^0(X_t,K_{\mathcal{F}_t})$ is not constant as $t\in \Delta$. 

Regarding foliations of Kodaira dimension 0,
in \cite{pereira05} Pereira has shown that if $\mathcal F$ is a foliation on a smooth surface $X$ such that $\kappa(\mathcal{F})=0$, then the smallest positive integer $k$
such that $h^0(X, kK_{\mathcal{F}})=1$ belongs to the set $\{1,2,3,4,5,6,8,10,12\}$.

In the proof of the theorem above, we extensively use McQuillan's results on the minimal model program for foliation on surfaces and the classification of the singularities of the canonical model of a foliation with pseudo-effective canonical divisor. 

As an immediate consequence of Theorem \ref{t_main} and Brunella's result above, we obtain the following:
\begin{cor}
Let $(X_t,\mathcal{F}_t)_{t\in\Delta}$ be a family of foliations with reduced singularities. Then, for any sufficiently large positive integer $m$, the dimension $h^0(X_t,\mathcal O_{X_t}(mK_{\mathcal F_t}))$ is constant for all $t\in \Delta$. 
\end{cor}

\section{Preliminary results}\label{s_preliminary}

We work over the field of complex numbers $\mathbb C$. 
We refer to \cite{KM98} for some of  the notations and basic results in birational geometry. 

A \textit{foliation} $\mathcal F$ on a $n$-dimensional smooth projective variety $X$ is given by a coherent subsheaf $T_{\mathcal F}$ of the tangent bundle $T_X$ of $X$ which is  closed under the Lie bracket and is such that  the quotient $T_X/T_{\mathcal F}$ is torsion free.  We denote by $\Omega_{\mathcal F}=T_{\mathcal F}^\ast$ the \textit{cotangent sheaf} of $\mathcal F$ and by $K_{\mathcal F}=c_1(\Omega_{\mathcal F})$ the \textit{canonical divisor} of $\mathcal F$. The \textit{singular locus} of $\mathcal F$, denoted by  ${\rm Sing} ~\mathcal F$, is defined as the set of points of $X$ on which $T_X/T_{\mathcal F}$ is not locally free. By Frobenius Theorem, around any point $p\in X$ outside the singular locus of $\mathcal F$, the germ of $T_{\mathcal F}$ coincides with the relative tangent bundle of a germ of a smooth fibration $X\supseteq U\to \mathbb C^q$, where $p\in U$. The \textit{dimension} of $\mathcal F$ is defined as $n-q$.

We now recall some basic facts about foliations over a surface (see \cite{brunella00} for more details).
Let $X$ be a smooth surface.
A foliation $\mathcal{F}$ of dimension 1  on $X$ is given by a short exact sequence
$$0\rightarrow T_{\mathcal{F}}\rightarrow T_X\rightarrow \mathcal{I}_Z N_{\mathcal{F}}\rightarrow 0$$
where $T_{\mathcal{F}}$ and $N_{\mathcal{F}}$ are line bundles and $\mathcal{I}_Z$ is an ideal sheaf supported on a finite set.
The line bundle $T_{\mathcal{F}}$ is the \textit{tangent bundle} of $\mathcal{F}$. The line bundle $N_{\mathcal{F}}$ is the \textit{normal bundle} of $\mathcal{F}$,
while its dual $N^{\ast}_{\mathcal{F}}$ is the \textit{conormal bundle}.
The support of $\mathcal O_X/\mathcal{I}_Z$ coincides with  ${\rm Sing}~\mathcal{F}$. 

Equivalently, a foliation on $X$ is the data of $\left\{(U_i,v_i)\right\}_{i\in I}$
where $\{U_i\}_{i\in I}$ is a covering of $X$, $v_i$ is a vector field on $U_i$ with only isolated zeroes
and there exist $g_{ij}\in \mathcal{O}^{\ast}(U_{i}\cap U_j)$ such 
that $$v_i\vert_{U_{i}\cap U_j}=g_{ij}v_j\vert_{U_{i}\cap U_j}
\qquad \text{for each }i,j\in I.$$
The cocycle $g_{ij}$ defines $K _{\mathcal{F}}$ as a line bundle.
A curve $C\subseteq X$ is said to be $\mathcal F$-\textit{invariant} 
if the inclusion $T_{\mathcal{F}}\vert_C\rightarrow T_X\vert_C$
factors through $T_C$.

If $p$ is a singular point of $\mathcal F$ and  $v$ is a vector field that defines $\mathcal F$ around $p$,
then  the eigenvalues of the linear part
$(Dv)(p)$ are defined up to multiplication by a non-zero constant. The point 
$p$ is  a \textit{reduced singularity} if at least one of the eigenvalues
of $(Dv)(p)$ is non-zero and their quotient is not a positive rational number.

Alternatively, a foliation $\mathcal F$ on a smooth surface $X$ can be locally defined by a holomorphic $1$-form $\omega$
with isolated zeroes (see \cite[pag. 19]{brunella00}). 

Let $\pi\colon\tilde{X}\rightarrow X$ be a proper birational morphism between smooth surfaces  and let $E$ be the exceptional curve.
Then the foliation $\mathcal{F}$ induces a foliation on $\tilde{X}\backslash E$ which
can be extended to a foliation on $\tilde{X}$ with isolated singularities.
We denote this foliation by $\pi^{\ast}\mathcal{F}$. 

\begin{teo}[Seidenberg] \cite[Chapt. 1, Thm. 1]{brunella00} \label{t_Seidenberg} Let $\mathcal F$ be a foliation on a smooth surface $X$. Then
for any $p\in{\rm Sing}(\mathcal{F})$,
there exists a sequence of blow-ups $\pi\colon \tilde{X}\rightarrow X$ over $p$
such that the foliation $\pi^{\ast}\mathcal{F}$ has only reduced singularities in a neighborhood of $\pi^{-1}(p)$.
\end{teo}

Given a foliation $\mathcal F$ on a smooth surface $X$, we define the \textit{Kodaira dimension of} $\mathcal F$ as the Kodaira dimension of $K_{\pi^{\ast}\mathcal F}$
 where $\pi$ is as in Theorem \ref{t_Seidenberg}
and we denote it by $\kappa(\mathcal F)$. It is easy to check that $\kappa(\mathcal F)$ does not depend on the resolution $\pi$. 
In particular, we say that $\mathcal F$ is of \textit{general type} if $\kappa(\mathcal F)=2$. Foliations of general type appeared several  times  in the literature, e.g. if $X=\mathbb P^2$, then the  foliations of general type on $X$ were studied by Pereira in 
\cite{pereira02}.

In this paper, we consider families of foliations defined over a smooth family of surfaces. 

\begin{defi}\cite[Def. 1]{brunella01}\label{hyp}
A \textit{family of foliations with reduced singularities} $(X_t,\mathcal{F}_t)_{t\in\Delta}$ 
(or a family of foliations, for short)
is the data of 
\begin{itemize}
\item a smooth morphism $\pi\colon\mathcal{X}\rightarrow\Delta$,
where $\mathcal{X}$ is a smooth complex
variety and $\Delta$ is the complex disc, whose fibres $X_t$ are projective surfaces, for all $t\in \Delta$;
\item a foliation $\mathcal{F}$ of dimension 1 on $\mathcal X$ such that
\begin{enumerate}
\item $\mathcal{F}$ is tangent to the fibres of $\pi$;
\item the singular set ${\rm Sing}~\mathcal{F}$ of $\mathcal{F}$
is of pure codimension 2 in $\mathcal{X}$ and cuts every fibre in a finite set; and
\item for any $t\in\Delta$, the foliation $\mathcal{F}_t=\mathcal{F}\vert_{X_t}$ is a foliation whose singularities 
${\rm Sing}~\mathcal{F}\vert_{X_t}={\rm Sing}~\mathcal{F}\cap X_t$
are reduced.
\end{enumerate}
\end{itemize}
\end{defi}
Note that \ref{hyp}(1) is needed to ensure that $\mathcal F_t$ is a foliation of dimension $1$ on $X_t$ for all $t\in \Delta$ and \ref{hyp}(2) is needed to ensure the existence of a canonical divisor of the foliation
$K_{\mathcal{F}}$ on $\mathcal{X}$ such that $$K_{\mathcal{F}}\vert_{X_t}=K_{\mathcal{F}_t}$$
for any $t\in\Delta$. 
Note also that invariance of the Kodaira dimension does not hold without hypothesis \ref{hyp}(3). Indeed, 
if the singularities of $\mathcal{F}_t$ are not reduced for all $t\in \Delta$, then 
invariance of plurigenera fails, as shown by the example in \cite[p. 114]{brunella01}. Furthermore, Example \ref{e} shows that the invariance of the Kodaira dimension does not hold without
the equality ${\rm Sing}~\mathcal{F}\vert_{X_t}={\rm Sing}~\mathcal{F}\cap X_t$: in the  example, we describe a family of foliations $\mathcal F_t$, induced by a foliation $\mathcal F$ on a smooth family of surfaces $\mathcal X\to \Delta$ such that $\mathcal F_t$ has reduced singularities for any $t\in \Delta$
and such that there exists a curve $C\subseteq X_{t_0}\cap {\rm Sing}~\mathcal{F}$, for some $t_0\in \Delta$.

On the other hand, under the assumptions above, we have: 
\begin{teo}\cite[Thm. 1]{brunella01}\label{bruinv}
Let $(X_t,\mathcal{F}_t)_{t\in\Delta}$ be a family of foliations on surfaces with reduced singularities.
Then the Kodaira dimension $\kappa(\mathcal{F}_t)$ does not depend on $t\in \Delta$.
\end{teo}
\medskip

Let $X$ be a  smooth surface. We will assume that  a curve $C\subseteq X$ is reduced and compact. 
We say that $C$ has \textit{normal crossing singularities} if locally, with respect to the Euclidean topology, around each point $p\in C$, the curve $C$ is a union of smooth curves meeting transversally. In particular a \textit{nodal} curve is an irreducible curve $C$ with normal crossing singularities. 

 For any curve $C\subseteq X$, the \textit{arithmetic Euler characteristic} of $C$ is given 
by $$\chi(C)=-K_X\cdot C-C^2.$$
Note that if $C$ is smooth, then  it coincides with the usual Euler characteristic.

Let $\mathcal F$ be a foliation on $X$ and let $p\in C$ be a point. If none of the components of $C$ is $\mathcal{F}$-invariant, 
we define the \textit{index of tangency} of $\mathcal{F}$ to $C$ at $p$
as follows.
Let $\{f=0\}$ be a local  equation of $C$ around $p$,
let $v$ be a local holomorphic vector field generating $\mathcal{F}$ around $p$.
Then, $${\rm tang}(\mathcal{F}, C, p)=\dim_{\mathbb{C}}\frac{\mathcal{O}_{X,p}}{<f,v(f)>}$$
where $v(f)$ is the Lie derivative of $f$ along $v$.
We have ${\rm tang}(\mathcal{F}, C, p)=0$ except on the finite subset of points of $C$
where $\mathcal{F}$ is not transverse to $C$.
Thus, we define  $${\rm tang}(\mathcal{F}, C)=\sum_{p\in C}{\rm tang}(\mathcal{F}, C, p).$$

\begin{prop}\cite[Chapt. 2, Prop. 2]{brunella00}\label{Cnotinv}
Let $\mathcal{F}$ be a foliation on a smooth surface $X$. 
Let $C$ be a curve on $X$ whose components are not $\mathcal F$-invariant.
Then
$$c_1(N_{\mathcal{F}})\cdot C=\chi(C)+{\rm tang}(\mathcal{F},C)\qquad\text{and}\qquad
K_{\mathcal{F}}\cdot C=-C\cdot C+{\rm tang}(\mathcal{F},C).$$
\end{prop}

\medskip

We now consider a curve  $C$  whose components are  all  $\mathcal{F}$-invariant.
If  $p\in C$ is a singular  point of $\mathcal F$, $\{f=0\}$  is a local equation for $C$ at $p$ and $\omega$ is a holomorphic 1-form that defines $\mathcal F$ around $p$, then we may write
 $$g\omega=h df+f\eta,$$
for some holomorphic 1-form $\eta$ and holomorphic functions $g,h$
such that $h$ and $f$ are coprime.
We define 
$$Z(\mathcal{F},C,p)=\text{vanishing order of }\left.\frac{h}{g}\right\vert_C\quad \text{at}\;\; p$$ and 
$$CS(\mathcal{F},C,p)=\text{residue of }-\left.\frac 1 h \eta\right\vert_C \quad \text{at }p.$$
If $p$ is a reduced singularity, then $Z(\mathcal{F},C,p)\geq 0$ by \cite{brunella99a}.
Let 
$$Z(\mathcal F,C)=\sum_{p\in\,{\rm Sing}(\mathcal{F})\cap C}Z(\mathcal{F},C,p)
\qquad \text {and}\qquad CS(\mathcal F,C)=\sum_{p\in\,{\rm Sing}(\mathcal{F})\cap C}CS(\mathcal{F},C,p).$$
\begin{prop}\label{Cinv}
Let $\mathcal{F}$ be a foliation on a smooth surface $X$ and  
let $C$ be a  curve on $X$ whose components are $\mathcal F$-invariant. 

Then
\begin{enumerate}
\item $c_1(N_{\mathcal{F}})\cdot C=C^2+Z(\mathcal{F},C)$;
\item $K_{\mathcal{F}}\cdot C=-\chi(C)+Z(\mathcal{F},C)$; and 
\item $C^2=CS(\mathcal{F},C)$.
\end{enumerate}
In addition, if $\mathcal F$ admits only reduced singularities, then $C$ has only normal crossing singularities.

\end{prop}
Formula (3) in Proposition \ref{Cinv} is usually referred as the \textit{Camacho-Sad formula}. 
\begin{proof}
(1) and (2) are implied by \cite[Chapt. 2, Prop. 3]{brunella00} 
and \cite{camachosad82} implies (3) (see also \cite[Chapt. 3, Thm. 2]{brunella00}).

Finally, if $\mathcal F$ admits only reduced singularities then  \cite[pag. 12]{brunella00} implies that $C$ has only normal crossing singularities.
\end{proof}

\medskip

\subsection{Minimal model program and Zariski decomposition for foliations}
In ~\cite{mcq08}, McQuillan has developed  a minimal model program for foliations on surfaces. We now recall some of the main results.

\begin{defi}\cite[Chapt. 5]{brunella00}
Let $\mathcal{F}$ be a foliation on a smooth surface $X$ and which admits only reduced singularities. 
We say that a curve $C$ in $X$ is $\mathcal{F}$-\textit{exceptional} if
\begin{enumerate}
\item $C$ is a smooth rational curve of self-intersection $-1$;
\item the contraction of $C$ to a point $p$ gives a new foliation $\overline {\mathcal F}$ such that  $p$ is either  a regular point or a reduced singular point for $\overline{\mathcal F}$.
\end{enumerate}

In particular, the foliation $\mathcal F$ is said to be \textit{relatively minimal} if 
$\mathcal F$ admits only reduced singularities and there are no $\mathcal F$-exceptional curves on $X$. 
\end{defi}

\medskip 

If $\mathcal F$ is a foliation on $X$ with reduced singularities, then there exists a birational morphism $X\to X'$ onto a smooth surface $X'$ such that the induced foliation $\mathcal F'$ on $X'$ is relatively minimal. Note that if $\pi\colon X\to \overline {X}$ is the contraction of an $\mathcal F$-exceptional curve $C$ onto a point $p$ and $\overline{\mathcal F}$ is the induced foliation on $\overline {X}$ then by \cite[pag. 72]{brunella00}, it follows that either $\mathcal F$ is regular at $p$ and $K_{\mathcal F}=\pi^{\ast}K_{\overline{\mathcal F}}+C$ or $\mathcal F$ is singular at $p$ and  $K_{\mathcal F}=\pi^{\ast}K_{\overline{\mathcal F}}$. 
Thus, we have that 
$$h^0(X,\mathcal O_{X}(mK_{\mathcal F}))=h^0(X',\mathcal O_{X'}(mK_{{\mathcal F'}}))$$ for all positive integers $m$. 

\medskip 

By the following result, $\mathcal{F}$-exceptional curves can be extended locally in a family:
\begin{lem}\label{Fecc}
Let $(X_t,\mathcal{F}_t)_{t\in\Delta}$ be a family of foliations and let $t_0\in \Delta$ be
such that $\mathcal{F}_{t_0}$ admits an $\mathcal{F}_{t_0}$-exceptional curve $E_{t_0}\subseteq X_{t_0}$.
Then there exists a neighborhood $t_0\in U\subseteq \Delta$
and a smooth hypersurface $E\subseteq \pi^{-1}(U)$ transverse to the fibres of $\pi$
such that $E_t=E\cap X_t$
is an $\mathcal{F}_t$-exceptional curve for any $t\in U$.

In particular, if $s\in \Delta$ then there exists a birational morphism $\nu\colon \mathcal X_U\to \mathcal X'_U$ over $U$ which defines a factorization
$$
\xymatrix{
\pi\colon \mathcal{X}_U \ar[r]^{\nu}&\mathcal{X}'_U \ar[r]^{\pi'} &U
}
$$
and such that the foliation $\mathcal{F}'$ induced on $\mathcal{X}'_U$ is relatively minimal on $X'_s$
(i.e. $\mathcal{F}'_s$ is relatively minimal) and 
the family $(X'_t,\mathcal{F}'_t)_{t\in U}$ is still a family of foliations with reduced singularities.
\end{lem}
\begin{proof}
The existence of $U$ and $E$ is guaranteed by \cite[Lemme 2]{brunella01}. 
Thus, after possibly shrinking $U$, 
there exists a birational morphism $\varepsilon \colon \mathcal X\to \mathcal Y$ which contracts $E$ and a smooth morphism $\pi'\colon \mathcal Y\to \Delta$ such that $\pi=\pi'\circ\varepsilon$ (see  \cite[Chapt. 1, Proof of Thm. 1.16]{fm94} for a similar argument).
For any $t\in U$, let $Y_t=\pi'^{-1}(t)$. Since $E|_{X_t}$ is $\mathcal F_t$-exceptional, the induced foliation on $Y_t$ admits reduced singularities for all $t\in U$. Thus, the claim follows after repeating the argument finitely many times. 
\end{proof}
\medskip 

We now consider a relatively minimal foliation $\mathcal F$ on a surface $X$ such that $K_{\mathcal F}$ is pseudo-effective. Then, we denote the \textit{Zariski decomposition} of $K_{\mathcal F}$ by 
$$K_{\mathcal F}=P+N$$
where $P$ is the positive part and $N$ is the negative part of $K_{\mathcal F}$. McQuillan shows that there exists a contraction $X\to X'$ onto a surface $X'$ with Kawamata log terminal singularities, which contracts all the curves contained in the support of $N$. More precisely, we say that a curve $C=\cup_{i=1}^r C_i$ in $X$ is an \textit{Hirzebruch-Jung string} if 
$C_1,\dots,C_r$ are smooth rational curves such that, for all $i,j\in \{1,\dots,r\}$, we have $C^2_i\le -2$, $C_i\cdot C_j =1$ if $|i-j|=1$ and $C_i\cdot C_j=0$ if $|i-j|>1$. Note that any Hirzebruch-Jung string on a smooth projective surface can be contracted onto a surface with cyclic quotient singularities (see \cite[Chapt. 3, Thm. 5.1 and Prop. 5.3]{bhpv04}).
We define: 

\begin{defi}\cite[Chapt. 8]{brunella00}
Given a foliation $\mathcal{F}$ on a surface $X$, we say that a curve $C$ is an $\mathcal{F}$-\textit{chain} if
\begin{enumerate}
\item $C=\cup_{i=1}^r C_i$ is an Hirzebruch-Jung string; 
\item each irreducible component $C_i$ of $C$ is $\mathcal{F}$-invariant;
\item the singularities of $\mathcal F$ along $C$ are all reduced; and
\item $Z(\mathcal{F},C_1)=1$ and $Z(\mathcal{F},C_i)=2$ for any $i\geq 2$.
\end{enumerate}
\end{defi}

\begin{teo}\cite[Thm. 2, Prop. III.1.2]{mcq08}, ~\cite[Chapt. 8, Thm. 1 and Addendum p. 109]{brunella00} \label{structneg}
Let $\mathcal{F}$ be a relatively minimal foliation on a smooth surface $X$, such that $K_{\mathcal F}$ is pseudo-effective. Let
$K_{\mathcal{F}}=P+N$ be the Zariski decomposition.

Then the support of $N$ is a disjoint union of maximal $\mathcal{F}$-chains and $\lfloor N\rfloor=0$. 
In particular, there exists a contraction $X\to X'$ onto a surface $X'$ with Kawamata log terminal singularities, which contracts all the curves in the support of $N$. 
\end{teo}

The following result is also due to McQuillan:
\begin{lem}\cite[Lemma IV.3.1]{mcq08}, \cite[Chapt. 8, Thm. 2,  and Chapt. 9, Thm. 2]{brunella00}\label{l_kod0}
Let $\mathcal F$ be a foliation on a smooth surface $X$. Assume that $\mathcal F$ admits only reduced singularities and that $\kappa(\mathcal F)=0$. Let $K_{\mathcal F}=P+N$ be the Zariski decomposition of $K_{\mathcal F}$. 

Then $P$ is a torsion divisor.  
\end{lem}

By a result of Brunella, the Zariski decomposition of the canonical divisor of a foliation
is well behaved in families:
\begin{prop}\cite[Prop. 1 and 2]{brunella01}\label{zardecfol}
Let $(X_t,\mathcal{F}_t)_{t\in\Delta}$ be a family of foliations of non-negative Kodaira dimension.
Then there exists an effective $\mathbb{Q}$-divisor $N$ on $\mathcal{X}$ such that $N$ does not contain any fibre of $\pi\colon \mathcal X\to \Delta$ and $N_t=N\vert_{X_t}$
is the negative part of the Zariski decomposition of $K_{\mathcal{F}_t}$ for any $t\in\Delta$.

In addition, if there exists $s\in \Delta$ such that $(X_s,\mathcal F_s)$ is relatively minimal, then there exists an open set $s\in U\subseteq \Delta$ such that the irreducible components $E_1,\dots,E_k$ of $N$ meet the surfaces $X_t$ transversally in distinct rational curve $E_1|_{X_t},\dots,E_k|_{X_t}$. 
\end{prop}

As a consequence, we have:

\begin{lem}\label{lcm}
Let $(X_t,\mathcal{F}_t)_{t\in\Delta}$ be a family of foliations of non-negative Kodaira dimension.
Let $N_t$  be the negative part of the Zariski decomposition of $K_{\mathcal{F}_t}$.

Then, the least common multiple of the denominators of $N_t$
does not depend on $t\in\Delta$.
\end{lem}
\begin{proof}
Let $m_t$ be the least common multiple of the denominators of $N_t$.
We will prove that $m_t$ is locally constant.
Fix $s\in \Delta$.
Let $U\subseteq \Delta$ be the neighborhood of $s$ and 
$$
\xymatrix{
\pi\colon \mathcal{X}_U \ar[r]^{\nu}&\mathcal{X}'_U \ar[r]^{\pi'} &U
}
$$
be the factorization of $\pi$ whose existence is guaranteed by Lemma \ref{Fecc}. Let $(X'_t,\mathcal F'_t)_{t\in U}$ be the induced family of foliations. 

By Proposition \ref{zardecfol}, modulo shrinking $U$,
there exists an effective $\mathbb{Q}$-divisor $N'$ on $\mathcal{X}'_U$
such that for any $t\in U$,
$N'_t=N|_{X'_t}$ is the negative part of the Zariski decomposition of $K_{\mathcal F_t}$ 
and the components of $N'_t$ meet the fibres of $\pi'$ transversally.
Thus, the least common multiple $m'_t$ of the denominators  of $N'_t$
does not depend on $t\in U$.

On the other hand, for any $t\in U$ the  foliation $\mathcal F'_t$ has only reduced singularities and in particular there exists an effective integral exceptional divisor $E_t$
such that $$K_{\mathcal{F}_t}=\nu_t^{\ast}K_{\mathcal{F}'_t}+E_t.$$
It follows that $N_t=\nu_t^{\ast}N'_t+E_t$ and the least common multiple of the denominators of $N_t$ coincides with $m'_t$. Thus, the claim follows. 
\end{proof}
%
%

Let $\mathcal{F}$ be a foliation on a smooth surface $X$.
Assume that $\mathcal{F}$ is of general type and let
$$K_{\mathcal{F}}=P+N$$
be the Zariski decomposition of $K_{\mathcal{F}}$.
We want to show that if $C$ is a  curve such that $P\cdot C=0$, then $C$ is $\mathcal{F}$-invariant. 
We assume first that $C$ is such that $K_{\mathcal{F}}\cdot C=0$ and assume by contradiction that $C$ is not $\mathcal{F}$-invariant.
Then, Proposition \ref{Cnotinv} implies 
$$0=K_{\mathcal{F}}\cdot C=-C^2+{\rm tang}(\mathcal{F},C)\geq-C^2$$
which is a contradiction because a big divisor cannot have intersection zero with a movable curve.
We now consider the general case. To this end, we consider a variation of ~\cite[Lemma III.1.1]{mcq08}:
\begin{prop}\label{lemmaneg}
Let $\mathcal{F}$ be a relatively minimal foliation on a smooth surface $X$ such that $K_{\mathcal F}$ is pseudo-effective. Let $\varepsilon \colon X\to Y$
be the contraction of all the components of the negative part of the Zariski decomposition of $K_{\mathcal F}$. 
Let $C$ be a  curve on $X$ which is not $\mathcal F$-invariant, let $\overline{C}$ be its image in $Y$ and let $\overline K=\varepsilon_*K_{\mathcal F}$.

Then $$(\overline K+\overline{C})\cdot \overline{C}\geq 0.$$
\end{prop}

We first prove the following: 

\begin{lem}\label{lemma}
Let $\varepsilon\colon X\to Y$ be a birational morphism between surfaces with only Kawamata log terminal singularities and assume that $\varepsilon$ contracts a chain of rational curves $F_1,\dots,F_k$. Let $L$ be a $\mathbb Q$-Cartier divisor on $X$ and let $C$ be  a curve on $X$ such that 
\begin{enumerate}
\item $(L+C)\cdot C\ge 0$;
\item $C\cdot F_i$ is a non-negative integer for any $i=1,\dots,k$; 
\item $F_1^2\le -1$, $F_i^2\le -2$ and $F_{i-1}\cdot F_{i}=1$ for any $i=2,\dots,k$; and 
\item $-1\le L\cdot F_1<0$ and $L\cdot F_i=0$ for any $i=2,\dots,k$.
\end{enumerate}

Let $\overline L=\varepsilon_*L$ and $\overline C=\varepsilon_* C$. Then $(\overline L+\overline C)\cdot \overline C\ge 0$. 
\end{lem}
\begin{proof}
We may write 
$$L+C=\varepsilon^{\ast}(\overline L+\overline C)-G$$
for some $\varepsilon$-exceptional divisor $G$. If $G\ge 0$ then the claim follows immediately. Therefore we may assume that $G$ is not effective and by the Negativity Lemma, there exists $i\in\{1,\dots,k\}$ such that $G\cdot F_i>0$ and the coefficient of $G$ along $F_i$ is negative. Thus, $(L+C)\cdot F_i<0$ and in particular we must have $i=1$ and $C\cdot F_1=0$. 
We may write $G=-\alpha F_1 + G_1$ for some $\alpha>0$ and some $\varepsilon$-exceptional divisor $G_1$ whose support does not contain $F_1$. 

There exists a morphism $f\colon X\to X_1$ which contracts only the curve $F_1$ and such that $\varepsilon$ factors through $f$. Let $L_1=f_*L$ and  $C_1=f_*C$.
 Since $L\cdot F_1<0$, we have $L=f^{\ast}L_1+\beta F_1$ for some $\beta >0$ and since $C\cdot F_1=0$ we have $C=f^{\ast}C_1$. In particular, 
$$(L_1+C_1)\cdot C_1=(L+C)\cdot C\ge 0.$$
Let $F'_i=f_*F_i$ for $i=2,\dots,k$. Then 
$C_1\cdot F'_i=C\cdot F_i$ for all $i=2,\dots,k$ and 
$L_1\cdot F'_i=0$ for any $i=3,\dots,k$. 
We may write $f^{\ast}F_2=\gamma F_1+ F_2$ for some $\gamma >0$. Then
$$0=F_1\cdot f^{\ast}F'_2=F_1\cdot (\gamma F_1+F_2)=-e_1\cdot\gamma+1,$$
where $e_1=-F^2_1\ge 1$. 
Thus, $\gamma=1/e_1$. 
We have,
$$F'_{i-1}\cdot F'_i=F_{i-1}\cdot F_i=1\qquad\text{and}\qquad {F'_i}^2=F_i^2$$
for any $i=3,\dots,k$. Moreover
$${F'_2}^2=F_2\cdot (\gamma F_1+F_2)=1/e_1 + F_2^2\le 1+ F_2^2 \le -1.$$
Finally,
$$L_1\cdot F'_2=L\cdot f^{\ast}F'_2=L\cdot (\gamma F_1+F_2)=1/e_1 L\cdot F_1\in [-1,0).$$

Thus, the claim follows by induction on $k$. 
\end{proof}

\begin{proof}[Proof of Proposition \ref{lemmaneg}]
We proceed by induction on the number of connected components of the exceptional locus of $\varepsilon$. 
By Theorem \ref{structneg}, any connected component of the exceptional locus of $\varepsilon$ is a maximal $\mathcal F$-chain $F_1,\dots,F_k$. 
Let $L=K_{\mathcal F}$. Then Proposition \ref{Cnotinv} implies that $(L+C)\cdot C\ge 0$. On the other hand, (2) of Proposition \ref{Cinv} implies that $L\cdot F_1=-1$ and $L\cdot F_i=0$ for any $i=2,\dots,k$. 
Thus, if $\varepsilon'\colon X\to Y'$ is the contraction of $F_1,\dots,F_k$, then Lemma \ref{lemma} implies that $(\varepsilon'_*L+\varepsilon'_*C)\cdot \varepsilon'_*C\ge 0$. 

Thus, the claim follows by induction, by proceeding as above for each connected component of the exceptional locus of $\varepsilon$. 
\end{proof}

\begin{prop}\label{zeroinv}
Let $\mathcal{F}$ be a relatively minimal foliation of general type on a smooth surface $X$.
Let $K_{\mathcal{F}}=P+N$ be the Zariski decomposition of $K_{\mathcal F}$ and
let $C$ be a curve such that $P\cdot C=0$.
Then $C$ is $\mathcal{F}$-invariant.
\end{prop}
\begin{proof}
Assume that $C$ is not $\mathcal F$-invariant. By Theorem \ref{structneg}, there exists a proper birational morphism $\varepsilon\colon X\rightarrow Y$ onto a normal surface $Y$ and whose exceptional locus coincides with the support of $N$.
Let $\overline K=\varepsilon_\ast K_{\mathcal F}$ and $\overline C=\varepsilon_*C$. 
The Negativity Lemma implies that $P=\varepsilon^{\ast}\overline K$. Thus,
$$\overline K\cdot \overline {C}=\varepsilon^*\overline K\cdot  C =P\cdot C =0.$$
Since $\overline K$ is big,
we have $\overline{C}^2<0$.
On the other hand,  Proposition \ref{lemmaneg} implies
$$\overline{C}^2=(\overline K+\overline{C})\cdot \overline{C}\geq 0$$
that is a contradiction.
Thus, $C$ is $\mathcal{F}$-invariant.
\end{proof}

The following Theorem is proved in \cite{mcq08}.
\begin{teo}\cite[Thm. 1 III.3.2, Remark III.2.2]{mcq08}\label{zeroloc}
Let $\mathcal{F}$ be a relatively minimal foliation of general type on a smooth surface $X$.
Let $K_{\mathcal{F}}=P+N$ be the Zariski decomposition of $K_{\mathcal F}$ and
let $Z$ be the union of all the curves $C$ such that  $P\cdot C=0$.
Then $Z$ is the union of:
\begin{enumerate}
\item the support of $N$;
\item disjoint chains of rational curves none of which is contained in the support of $N$;
\item cycles $\Gamma$ of rational curves such that ${\rm Sing}(\mathcal{F})\cap \Gamma$ coincides with the singular locus of
$\Gamma$; and
\item single rational nodal curves $\Gamma$ such that ${\rm Sing}(\mathcal{F})\cap \Gamma$ coincides with the singular locus of $\Gamma$.
\end{enumerate}
Moreover, a chain $\mathcal{C}$ of type $(2)$ is either disjoint from ${\rm Supp}N$
or there exist exactly two connected components of ${\rm Supp} N$, each of which  consists of a smooth rational curve $E_i$ of self-intersection $-2$, with $i=1,2$, and such that both $E_1$ and $E_2$ meet $\mathcal C$  transversally along the same tail $C$ of $\mathcal C$ on the points $p_1$ and $p_2$, so that the intersection is transverse and $N\vert_C=\frac 12 p_1+\frac 12 p_2$.
\end{teo}

\begin{oss}\label{noint}
Under the same assumptions as in Theorem \ref{zeroloc},  it follows by Proposition \ref{zeroinv}  that 
any curve $C$ such that $P\cdot C=0$ is $\mathcal{F}$-invariant.
Since any point lying in the intersection of two $\mathcal{F}$-invariant
curves is a singular point, 
the cycles of rational curves and the nodal curve appearing in
Theorem \ref{zeroloc}
do not meet any other $\mathcal{F}$-invariant curve.
In particular they do not meet any component of ${\rm Supp}N$.
\end{oss}

\begin{defi}\label{d_egl}
A connected component of type (3) and (4) of Theorem \ref{zeroloc}
is called an \textit{elliptic Gorenstein leaf} (e.g.l. for short).
\end{defi}

The following theorem due to McQuillan will play a key role in the proof of our main results:

\begin{teo}\cite[Thm. IV.2.2]{mcq08}\label{contrEGL}
Let $\mathcal{F}$ be a relatively minimal foliation on a smooth surface $X$.
Let $\Gamma$ be an elliptic Gorenstein leaf.
Then $K_{\mathcal{F}}\vert_{\Gamma}$ is not a torsion divisor.
\end{teo}

%

\medskip

\subsection{Foliations and fibrations}
When a surface is endowed with a fibration,
the study of foliations on the variety becomes simpler.
In particular two types of foliations play a key role in the case of Kodaira dimension one:
foliations induced by fibrations and foliations transverse to fibrations.

\subsubsection{Elliptic fibrations}
We first recall some of the basic notions for the canonical bundle formula for an elliptic fibration (e.g. see  \cite[pag. 236]{ambro04} for more details). Let $X$ be a smooth surface and let $f\colon X\to C$ be an elliptic fibration onto a curve $C$. Let
$f'\colon X'\rightarrow C$ be the relatively minimal elliptic fibration associated to $f$, obtained by blowing-down any possible sequence of vertical $(-1)$-curves. Thus, we obtain a diagram:  
$$
\xymatrix{
X \ar[r]^{\varepsilon}\ar[rd]_f& X' \ar[d]^{f'}\\
 &C.
}
$$
The \textit{discriminant} of $f$ is defined by 
\begin{equation}\label{e_disc}
B_C=\sum_{p\in C} (1-\gamma_p)p
\end{equation}
where, for any $p\in C$, $\gamma_p$ denotes the \textit{log canonical threshold} of  $X'$ with respect to $f'^*p$:
 $$\gamma_p=\sup\{t\in \mathbb{R}_{>0}\mid (X',t f'^{\ast}p)\text{ is log canonical}\}.$$
Then $K_{X'/C}=f'^{\ast}(M_C+B_C)$
where $M_C$ is a $\mathbb Q$-divisor on $C$ which denotes the \textit{moduli part} in the canonical bundle formula  of $f$ \cite{kodaira64}. In particular, $\deg M_C\ge 0$ and the equality holds if and only if $f$ is isotrivial. 
 
We may write, 
\begin{equation}\label{e_disc2}
B_C=B'_C+\sum \frac{m_p-1}{m_p}p
\end{equation}
where the sum runs over the points $p$ such that
the fibre $F$ of $f'$ over $p$ is a multiple fibre of multiplicity $m_p$,
that is, $f'^{\ast}p=m_p F_{red}$, where $F_{red}$ is the reduced divisor associated to $f'^{\ast}p$.
Then $$L=M_C+B'_C$$ is an integral divisor on $C$.
Thus, if we denote by $\{M_C\}$ the fractional part of $M_C$, then
\begin{eqnarray}\label{fracpart}
{\rm Supp}\{M_C\}={\rm Supp}B'_C.
\end{eqnarray}
Furthermore, by \cite{kodaira64}, it follows that 
\begin{equation}\label{e_12cartier}
12M_{C}\text{ is  Cartier and }|12M_C| \text{ is base point free}. 
\end{equation}

We will often use Kodaira's classification of the singular fibres of an elliptic fibration. In particular, 
if $p\in C$ is such that the fibre $f^{\ast}p$ is singular and $b_p=1-\gamma_p$ is the coefficient of $B_C$ along $p$, then the fibre $f^{-1}(p)$ is of one of the following types: 
$$ mI_b, I^{\ast}_b, II, II^{\ast}, III, III^{\ast}, IV, IV^{\ast}
$$
and the corresponding  values of $b_p$ are:
$$1-\frac 1 m, \frac 1 2, \frac 1 6, \frac 5 6, \frac 1 4,\frac 3 4, \frac 1 3, \frac 2 3.$$
\begin{oss}\label{r_lct}
Let $X$ be a smooth surface and let $f\colon X\to C$ be an elliptic fibration onto a curve $C$ and assume that all the fibres of $f$ have support with only normal crossing singularities. Let $\varepsilon\colon X\to X'$ be the relative minimal elliptic surface associated to $f$ and let $f'\colon X'\to C$ be the induced fibration.
We may write $K_X=\varepsilon^*K_{X'}+E$ where $E=\sum a_D D\ge 0$ is an $\varepsilon$-exceptional divisor. Let $p\in C$. We may write 
$$f^*p=\sum_{D\subseteq f^{-1}(p)} l_DD.$$ Then  
 the log canonical threshold $\gamma_p$
is computed on $X$ and it is equal to
$$\gamma_p=\min\left\{\left.\frac{1+a_D}{l_D}\right\vert D\subseteq f^{-1}(p)\right\}.$$
We say that a prime divisor $D\subseteq f^{-1}(p)$ {\em computes the log canonical threshold} at $p$ if $\frac{1+a_D}{l_D}=\gamma_p$.  
\end{oss}

\medskip

\subsubsection{Foliations induced by fibrations}
Let $X$ be a smooth surface and let $f\colon X\rightarrow C$ be a fibration onto a smooth curve.
The fibration induces a foliation $\mathcal{F}$ whose leaves are contained in the fibres of $f$.
The canonical divisor of $\mathcal{F}$ is (cf. \cite[Chapt. 2, Section 3]{brunella00})
\begin{eqnarray}\label{folfibr}
K_{\mathcal{F}}&=&K_{X/C}+\sum(1-l_D)D
\end{eqnarray}
where  the sum runs over all the irreducible curves $D$ contracted by $f$ and $l_D$ denotes the ramification order of $f$ along $D$, i.e. for any $p\in C$ we have $$f^{\ast}p=\sum_{D\subseteq f^{-1}(p)}l_D D.$$
\begin{oss}\label{r_redsnc}
Let  $\mathcal F$ be a foliation on a smooth surface $X$ induced by a fibration $f\colon X\to C$ onto a curve $C$. Then
 $\mathcal F$ admits only  reduced singularities
if and only if all the fibres of $f$ have support with only normal crossing singularities. 
Indeed if $\mathcal F$ is reduced then by Proposition \ref{Cinv} all the $\mathcal F$-invariant curves have support with normal crossing singularities. The other direction follows from an explicit computation.
\end{oss}

If a foliation $\mathcal{F}$ is induced by an elliptic fibration, then
we can give a  precise description of the Zariski decomposition of $K_{\mathcal{F}}$.
\begin{lem}\label{cbf} Let $X$ be a smooth surface and let $f\colon X\rightarrow C$ be a fibration onto a smooth curve. 
Let $\mathcal{F}$ be the  foliation induced by $f$. 
\begin{enumerate}
\item For any $p\in C$, there exists a neighborhood $\mathcal U$ of $f^{-1}(p)$ such that 
$$K_{\mathcal F}|_{\mathcal U}\sim_{\mathbb Q}(K_X+F_{red})|_{\mathcal U}$$
where $F_{red}$ denotes the reduced divisor associated to $f^*p$.  
\item 
Moreover, if $f$ is an elliptic fibration,  $\mathcal{F}$ admits only reduced singularities,
 $K_{\mathcal{F}}$ is pseudo-effective and $K_{\mathcal{F}}=P+N$
is the Zariski decomposition of $K_{\mathcal{F}}$, then
$$P=f^{\ast}M_C$$ where $M_C$ is the moduli part in the canonical bundle formula.
\end{enumerate}
\end{lem}

\begin{proof}
Let $p\in C$ and let $\mathcal U$ be a sufficiently small neighborhood of $f^{-1}(p)$. If $l_D$ denotes the ramification order of $f$ along $D$ for all $D\subseteq f^{-1}(p)$,  then   \eqref{folfibr} implies
$$K_{\mathcal F}|_{\mathcal U}\sim_{\mathbb Q}( K_X + \sum_{D\subseteq f^{-1}(p)} (1-l_D)D)|_{\mathcal U}=(K_X+F_{red})|_{\mathcal U}$$
and $(1)$ follows. 

Let us assume now that $f$ is as in (2). Let $\varepsilon\colon X\to X'$ be the relative minimal elliptic surface associated to $f$ and let $f'\colon X'\to C$ be the induced fibration. 

Then $K_{X'/C}=f'^{\ast}(M_C+B_C)$
where $M_C$ is the moduli part in the canonical bundle formula and 
$B_C$ is the discriminant.
Since $\mathcal F$ has only reduced singularities,  Proposition \ref{Cinv}  implies that the fibres of $f$ have support with only normal crossing singularities. Thus, as in Remark \ref{r_lct}, the log canonical threshold $\gamma_p$ is computed on $X$ so that 
$$\gamma_p=\min\left\{\left.\frac{1+a_D}{l_D}\right\vert D\subseteq f^{-1}(p)\right\}$$
where $E\ge 0$ is an $\varepsilon$-exceptional divisor such that  
$K_X=\varepsilon^{\ast}K_{X'}+E$ and $E=\sum a_D D$.

Then \eqref{folfibr} implies that 
$$
\begin{aligned}
K_{\mathcal{F}}&=K_{X/C} + \sum(1-l_D)D\\
&=f^{\ast}(M_C+B_C)+ E+ \sum(1-l_D)D.
\end{aligned}
$$

We claim that  $$\Psi=f^{\ast}(B_C)+ E+ \sum(1-l_D)D$$
is an effective divisor whose support does not contain any fibre. Indeed, for any $p\in C$ and for any prime divisor $D\subseteq f^{-1}(p)$,  
the coefficient of $\Psi$ along $D$ is
 $$1+a_D-l_D\gamma_p\ge 0$$
 and the equality holds for any $D$ which computes the log canonical threshold (cf. Remark \ref{r_lct}). 
 It follows that $\Psi\ge 0$ and, for any $p\in C$, there exists a prime divisor $D\subseteq f^{-1}(p)$ which is not contained in the support of $\Psi$. Thus, the claim follows. 
 
Since ${\rm Supp} \Psi$ does not contain any fibre, we have that $\Psi=N$ and $P=f^{\ast}M_C$. Thus, $(2)$ follows. 
\end{proof}

\medskip

\subsubsection{Foliations transverse to a fibration}
Let  $X$ be a smooth surface and let $f\colon X\rightarrow C$ be a fibration onto a curve $C$. Let $\mathcal{F}$ be a foliation on $X$
which is \textit{transverse} to $f$, that is, such that the general fibre $F$ of $f$ is not $\mathcal{F}$-invariant
and $K_{\mathcal{F}}\cdot F=0$.
Thus, there exists an effective divisor $D_{\tan}$ (cf. \cite[p. 573, Lemme 4]{brunella97}), whose support is contained in  the set of  $\mathcal{F}$-invariant curves 
contained in the fibres of $f$ and 
 such that
\begin{eqnarray}\label{cantrasv}
K_{\mathcal{F}}=f^{\ast}K_C + D_{\tan}+\sum (l_D-1)D
\end{eqnarray}
where  the sum runs over all the irreducible curves $D$ contracted by $f$ and $l_D$ denotes the ramification order of $f$ along $D$, i.e. for any $p\in C$ we have $$f^{\ast}p=\sum_{D\subseteq f^{-1}(p)}l_D D.$$
Let $D_f=\sum (l_D-1)D$.
Since $D_{\tan}+D_f$ is contained in fibres of $f$,
its Zariski decomposition is 
\begin{eqnarray}\label{Dtan}
D_{\tan}+D_f=&f^{\ast}\theta+\bar{N}
\end{eqnarray}
where $\bar{N}$ is the negative part of the Zariski decomposition and  $\theta$ is the largest effective $\mathbb Q$-divisor such that $D_{\tan}+D_f-f^{\ast} \theta$ is effective.
In particular, 
$$K_{\mathcal F}=f^{\ast}(K_C+\theta)+{\overline N}.$$
Since the support of $\overline N$  is contained in fibres of $f$ but it does not contain any of its fibres, it follows that if $K_{\mathcal{F}}$ is pseudo-effective and $K_{\mathcal{F}}=P+N$
is its Zariski decomposition then $K_C+\theta$ has non-negative degree and  $\bar{N}=N$.

\begin{lem}\label{D_tan} With the notation introduced above, 
let $\theta=\sum_{q\in C} \theta_q q$ and let $p\in C$ be such that $\theta_p\in\mathbb{Q}\backslash\mathbb{Z}$.
Then the support of $\bar{N}$ contains all the components $D$ of $f^{\ast}p$ that are reduced (i.e. such that $l_D=1$).
\end{lem}
\begin{proof}
We may write $$f^{\ast}p=\sum l_D D\qquad D_{\tan}=\sum a_D D\qquad\text{and}\qquad \bar{N}=\sum c_D D.$$
Let $\bar{D}$ be a reduced component of $f^{\ast}p$. 
Then \eqref{Dtan} implies 
$$a_{\bar{D}}=\theta_p+c_{\bar{D}}.$$
Since $a_{\bar{D}}$ is an integer and $\theta_p$ is not, it follows that $c_{\bar{D}}\neq 0$, as claimed. 
\end{proof}

\begin{defi}\label{d_riccati}
A foliation $\mathcal{F}$ on a smooth surface $X$ and which is transverse to a fibration  $f\colon X\rightarrow C$
whose general fibre is a rational curve is called a \textit{Riccati foliation}.
\end{defi}

\begin{defi}\label{d_turbulent}
A foliation $\mathcal{F}$ on a smooth surface $X$ and which is  transverse to an  elliptic fibration $f\colon X\rightarrow C$
 is called a \textit{turbulent foliation}.
\end{defi}

\begin{oss}\label{r_turbulent}
If $\mathcal F$ is a turbulent foliation associated to an elliptic fibration $f\colon X\rightarrow C$, then $f$ is isotrivial
\cite[pag. 64]{brunella00}. In particular, it follows that $f$ does not admit fibres of type $I_b$ and $I_b^{\ast}$ for $b\ge 1$
\cite[pag. 68]{brunella00}.
\end{oss}

Similarly to (1) of Lemma \ref{cbf}, we have:

\begin{prop}\label{prop6ii}\cite[pag. 69-70]{brunella00}.
Let $X$ be a smooth surface and let  $\mathcal{F}$ be a  turbulent foliation on $X$ which is transverse to the elliptic fibration $f\colon X\rightarrow C$. Assume that $\mathcal F$ admits only reduced singularities. Let $p\in C$ and let $F_{red}$ be the reduced divisor associated to $f^{\ast}p$. 

Then $F_{red}$ has normal crossing singularities and there exists a neighborhood $\mathcal U$ of $f^{-1}(p)$ such that 
$$K_{\mathcal F}|_{\mathcal U}\sim_{\mathbb Q}(K_X+F_{red})|_{\mathcal U}.$$ 
\end{prop}

\begin{cor}\label{c_turb} Let $X$ be a smooth surface which admits an elliptic fibration $f\colon X\to C$ and let $\varepsilon \colon X\to X'$ be the relative minimal fibration associated to $f$ with induced fibration $f'\colon X'\to C$. 
Let $\mathcal F$ be either a turbulent foliation on $X$ which is transverse to $f$ or the foliation induced by $f$. 
Let $B_C$ be the discriminant of $f$ and let $B'_C$ be the $\mathbb Q$-divisor on $C$ defined in \eqref{e_disc2}. Assume that $\mathcal F$ admits only reduced singularities and that $K_{\mathcal F}$ is pseudo-effective. Let $K_{\mathcal F}=P+N$ be the Zariski decomposition of $K_{\mathcal F}$ and assume that $\lfloor N\rfloor =0$. 

Then 
\begin{enumerate}
\item  $f({\rm Supp }N)={\rm Supp }B'_C$ and in particular $N= 0$ if and only if $B'_C=0$; 
\item  if $p\in {\rm Supp } B'_C$ then the coefficient of $B_C$ at $p$ coincides with the coefficient of $N$ along any reduced component $E$ of $f^{\ast}p$ which is not contained in the exceptional locus of $\varepsilon$;
\item if $p\in {\rm Supp } B'_C$ and $D\subseteq f^{-1}(p)$ is a prime divisor, then $D$ is contained in the support of $N$ if and only if it does not compute the log canonical threshold at $p$ (cf. Remark \ref{r_lct}); and
\item if $\mathcal F$ is turbulent and  $f^{\ast} p$ is a multiple fibre for some $p\in C$, then the reduced divisor associated to $F$ is the union of a smooth curve of genus one and trees of $\mathcal F$-exceptional curves. 
\end{enumerate}
\end{cor}
\begin{proof}
 We have $K_X=\varepsilon^{\ast}K_{X'}+E$ for some $\varepsilon$-exceptional divisor $E=\sum a_DD\ge 0$.
By $(1)$ of Lemma \ref{cbf} and Proposition \ref{prop6ii}, it follows that for all $p\in C$, there exists a neighborhood $\mathcal U$ of $f^{-1}(p)$ such that 
$$K_{\mathcal F}|_{\mathcal U}\sim_{\mathbb Q} (K_X+F_{red})|_{\mathcal U}$$
where $F_{red}$ denotes the reduced divisor associated to $f^*p=\sum_{D\subseteq f^{-1}(p)}l_D D$.  
Since $\varepsilon^*K_{X'}|_{\mathcal U}\sim_{\mathbb Q} 0$, it follows that if $\gamma_p$ is the log canonical threshold of $X'$ with respect to $f'^{\ast}p$, then
$$K_{\mathcal F}|_{\mathcal U}\sim_{\mathbb Q} 
(\sum_{D\subseteq f^{-1}(p)} (a_D+1)D)|_{\mathcal U}\sim_{\mathbb Q} 
(\sum_{D\subseteq f^{-1}(p)} (a_D+1-\gamma_p l_D)D)|_{\mathcal U}.
$$
On the other hand, we have that since $\mathcal F$ admits only reduced singularities, $F_{red}$ has  normal crossing singularities and as in Remark \ref{r_lct}, we have
$$\gamma_p=\min\left\{\left.\frac{1+a_D}{l_D}\right\vert D\subseteq f^{-1}(p)\right\}.$$
By (2) of Lemma \ref{cbf} and \eqref{cantrasv} it follows that the Iitaka fibration of $K_{\mathcal F}$ factors through $f$. Thus, $P|_{\mathcal U}\sim_{\mathbb Q}0$ and 
$$N|_{\mathcal U}=(\sum_{D\subseteq f^{-1}(p)} (a_D+1-\gamma_p l_D)D)|_{\mathcal U}.$$
 
 If $p\notin {\rm Supp} B'_C$ then $\gamma_p=\frac 1 {m_p}$, where $m_p$ denotes the multiplicity of $f'^{\ast}p$. Thus,  $\gamma_pl_D$ is a positive integer for each $D\subseteq f^{-1}(p)$.   In particular, since by assumption $\lfloor N\rfloor =0$, we have that    $N|_{\mathcal U}= 0$. 
 
On the other hand, if $p\in {\rm Supp } B'_C$ then by the classification of the singular fibres of an elliptic fibration, there exists a component $D_0$ of $F_{red}$ such that $l_{D_0}=1$ and $a_{D_0}=0$. Thus $(1)$ and $(3)$ follow. 
If  $E\subseteq f^{-1}(p)$ is a reduced component of $f^*p$ which is not contained in the exceptional locus of $\varepsilon$, we 
have that $a_{E}=0$ and $l_E=1$.  
In particular, the coefficient of $N$ along $E$ is $1-\gamma_p$. Thus $(2)$ follows.

Assume now that $\mathcal F$ is turbulent and $f^{\ast}p$ is a multiple fibre for some $p\in C$ and let $F_{red}$ be the reduced divisor associated to $F$. Suppose  that $F_{red}$ is not irreducible. By Remark \ref{r_turbulent}, $f^{-1}(p)$ is not of type $I_b$ or $I_b^*$ for $b\ge 1$. 
Then 
by the classification of the singular fibres of an elliptic fibration, 
$f^{\ast}p$ contains a smooth curve of genus 1 in its support and
there exists a $(-1)$-curve $E_0$ contained in $f^{-1}(p)$  which meets $F_{red}-E_0$ transversally in either one or two points. 
 Proposition \ref{prop6ii} implies 
$$K_{\mathcal F}\cdot E_0=(K_X+F_{red})\cdot E_0= -2+(F_{red}-E_0)\cdot E_0\le 0.$$ 
In particular, Proposition \ref{Cnotinv} implies that $E_0$ is $\mathcal F$-invariant, as otherwise
$$K_{\mathcal F}\cdot E_0=-E^2_0+{\rm tang}(\mathcal{F},E_0)\ge 1.$$

 Thus, Proposition \ref{Cinv} implies that  $Z(\mathcal F,E_0)=(F_{red}-E_0)\cdot E_0\le 2$ and $E_0$ is $\mathcal F$-exceptional \cite[pag. 72]{brunella00}. Let $g\colon X\to Y$ be the contraction of $E_0$ and let $\mathcal F'$ be the induced foliation on $Y$. Then, after replacing $X$ by $Y$ and $\mathcal F$ by $\mathcal F'$ and repeating the same argument as above finitely many times,  we may assume that there exists a $(-1)$-curve $E_0$ contained in $f^{-1}(p)$ such that $(F_{red}-E_0)\cdot E_0=1$ and in particular $K_{\mathcal F}=g^*K_{\mathcal F'}+E_0$. Thus, $E_0$ is contained in the support of $\lfloor N\rfloor$, a contradiction. Thus, $F_{red}$ is a smooth curve of genus one and $(4)$ follows. 
 \end{proof}

\medskip

\subsubsection{Isotrivial fibrations}
Given a smooth surface $X$, we consider an isotrivial fibration
$$\psi\colon X\rightarrow D$$ over a curve $D$ and whose general fibre $F$ has genus 
greater than one. Then there exists a curve $G$ and a finite group $\Gamma$ acting  on $G$ and $F$ such that, if we consider the diagonal action of $\Gamma$ on $G\times F$, then $D=G/\Gamma$,  
$X$ is birational to $ (G\times F)/\Gamma$ and the 
induced diagram 
$$
\xymatrix{
X\ar[rd]_{\psi}\ar@{-->}[r] &(G\times F)/\Gamma\ar[d]^{h}\\
&D=G/\Gamma.
}
$$
is commutative.  \cite[Prop 2.2]{serrano96}\label{h1} implies that 
\begin{equation}\label{e_h1}
h^1(X,\mathcal{O}_X)=g(G/\Gamma)+g(F/\Gamma).
 \end{equation}

Let $\varepsilon \colon Y\to (G\times F)/\Gamma$ be the minimal resolution and let $f\colon Y\to D$ be the induced morphism. Then the exceptional locus of $\varepsilon$ is a disjoint union of  Hirzebruch-Jung strings, each of which meets the strict transform of a fibre of $\psi$ transversally in one point \cite[Theorem 2.1]{serrano96}.

\begin{oss}\label{impl2} With the notation introduced above, 
we claim that $f\colon Y\to D$ is the minimal fibration so that each fibre has support with only normal crossing singularities, which means that if $f'\colon Y'\to D$ is any fibration from a smooth surface $Y'$ which is  birational to $f\colon Y\to D$, i.e. there exists a birational map $\eta\colon Y'\dashrightarrow Y$ such that $f\circ \eta = f'$,  and such that each fibre of $f'$ has support with only normal crossing singularities, then there exists a proper birational morphism $Y'\to Y$ which defines a factorization of $f'$. 

Indeed, let $q\colon W\to Y$ and $q'\colon W\to Y'$  be proper  birational morphisms from a smooth surface $W$ which resolve the 
indeterminacy of $\eta$, i.e. $\eta \circ  q' = q$. It is enough to show that the exceptional locus of $q'$ is contained in the exceptional locus of $\varepsilon\circ q$, which implies that there exists a proper birational morphism $Y'\to (G\times F)/\Gamma$ and the claim follows from the fact that $Y$ is the minimal resolution of $(G\times F)/\Gamma$. 

Let $h\colon (G\times F)/\Gamma\to D$ be the induced moprhism and let us assume by contradiction that the exceptional locus of $q'$ is not contained in the exceptional locus of $\varepsilon\circ q$. It follows that there exists $p\in D$ such that if $E$ is the strict transform in $Y'$
of the support $E'$ of the fibre $h^{-1}(p)$ then $E$ is $q'$-exceptional.
In particular $E$ is a rational curve. 
Since $F$ has genus greater than one, there are at least three singular points of $(G\times F)/\Gamma$ along $E'$ which coincide with the branch points of the induced finite morphism $F\to E'$. Thus, there exist at least three Hirzebruch-Jung strings inside $f^{-1}(p)$ intersecting $E$ and in particular the fibre $f'^{-1}(p)$ does not have support with only normal crossing, a contradiction. 
Thus, the claim follows. 
\end{oss}

Let us assume now that any fibre of  $\psi\colon X\to D$ has  support with only normal crossing singularities.  Remark \ref{impl2} implies that  $\psi$ factors through  a proper birational morphism $X\to (G\times F)/\Gamma$.
Let $C=F/\Gamma$. Then the induced morphism $\varphi\colon X\to C$  is a fibration which is called {\em transverse} to $\psi$. 
Note that $\varphi$ is also an isotrivial fibration with general fibre isomorphic to $G$. From now on, we assume that the genus of $G$ is greater than one. 

As above, we denote by $\varepsilon \colon Y\to (G\times F)/\Gamma$ the minimal resolution and by $\psi'\colon Y\to D$ and $\varphi'\colon Y\to C$  the induced morphisms. By Remark \ref{impl2}, there exists a proper birational morphism $\nu\colon X\to Y$. 

For any $p\in D$ and for any prime divisor $E\subseteq \psi^{-1}(p)$,  we denote by $l_E$ the ramification index of $\psi$ along $E$, so that  $\psi^{\ast}p=\sum l_E E$. Thus, we define $D_\psi=\sum (l_E-1)E$. Similarly, we define $D_{\varphi}$, $D_{\psi'}$ and $D_{\varphi'}$. Note that $D_{\varphi'}=\nu_*D_{\varphi}$ and $D_{\psi'}=\nu_*D_\psi$. 

\begin{teo}\cite[Thm. 2.1(i) and Thm. 4.1]{serrano96}\label{serr41}
With the notation introduced above, we have that the support of any fibre of $\varphi$ and $\psi$ has support with only normal crossing singularities.

Further, if $Z$ is the reduced divisor on $Y$ whose support coincides with the exceptional locus of $\varepsilon$, we have 
$$K_Y=\varphi'^{\ast}K_{C}+D_{\varphi'}+ \psi'^{\ast}K_{D}+D_{\psi'}+Z.$$
\end{teo}
\begin{lem} \label{iitaka} With the notation introduced above, 
we have that
the morphism $\varphi'$ is the Iitaka fibration of $K_{Y/D}-D_{\psi'}$
and  the morphism $\psi'$ is the Iitaka fibration of $K_{Y/C}-D_{\varphi'}$.
\end{lem}
\begin{proof}
 \cite[Prop. 3.1 and Prop. 5.1]{serrano96} implies that $$\kappa(K_{Y/D}- D_{\psi'})=1.$$
Let $G'$ be the general fibre of $\varphi'$. Then,  if $Z$ is the reduced divisor on $Y$ whose support coincides with the exceptional locus of $\varepsilon$, Theorem \ref{serr41} implies 
$$G'\cdot( K_{Y/D}- D_{\psi'})= G'\cdot(\varphi'^{\ast}K_{C}+ D_{\varphi'}+Z)=0.$$
Thus, it follows that $\varphi'$ is the Iitaka fibration of $K_{Y/D}-D_{\psi'}$. Similarly,  $\psi'$ is the Iitaka fibration of $K_{Y/C}-D_{\varphi'}$.
\end{proof}

Thus, given an isotrivial fibration $\psi\colon X\to D$ as above and the associated transverse fibration $\varphi\colon X\to C$, there are two foliations on $X$ which are naturally associated to $\psi$. Indeed,  we denote by $\mathcal F$  the foliation induced by $\psi$ and by $\mathcal G$ the foliation induced by $\varphi$. 
By \eqref{folfibr}, we have
\begin{equation}\label{e_kiso}
K_{\mathcal F}=K_{X/D}-D_\psi\qquad\text{and}\qquad K_{\mathcal G}=K_{X/C}-D_\varphi.
\end{equation}
In particular, Lemma \ref{iitaka} and Proposition \ref{Cnotinv} imply that ${\rm tang}(\mathcal F,G)=0$ and therefore 
  $\mathcal F$ is transverse to $\varphi$. Similarly, $\mathcal G$ is transverse to $\psi$.
By  Remark \ref{r_redsnc} and Theorem \ref{serr41}, $\mathcal F$ and $\mathcal G$ have reduced singularities.


\begin{lem}\label{overpoint}  With the notation introduced above, we have that 
\begin{enumerate}

\item  $\mathcal{F}$ and $\mathcal G$ are relatively minimal if and only if
$X$ is the minimal resolution of $(G\times F)/\Gamma$;
\item the support of the negative part in the Zariski decomposition of $K_{\mathcal F}$ coincides with the exceptional locus of the morphism $X\to (G\times F)/\Gamma$; and 
\item if $p\in D$ then there exists a unique component $E$ of the fibre $\psi^*p$ 
which is not $\mathcal{F}$-invariant and all the other components are contained in the exceptional locus of the morphism $X\to (G\times F)/\Gamma$.
\end{enumerate}
\end{lem}
\begin{proof}

 Remark \ref{impl2} implies that the minimal resolution of $(G\times F)/\Gamma$
 is the minimal fibration so that each fibre has support with only normal crossing singularities. Thus, 
 by Remark \ref{r_redsnc}, $(1)$ follows. 

We now show $(2)$. We first assume that $\mathcal F$ is relatively minimal. By $(1)$, we have that $\varepsilon\colon X\to (G\times F)/\Gamma$ is the minimal resolution.
Let $C_1,\dots,C_m$  be a maximal Hirzebruch-Jung string contained in the exceptional locus of $\varepsilon$. Then, only one tail, say $C_m$ intersects the rest of the fibre and,  in particular, we have that $C_1,\dots,C_m$ are $\mathcal F$-invariant and such that    $Z(\mathcal F,C_1)=1$ and $Z(\mathcal F,C_i)=2$ if $i=2,\dots,m$. Proposition \ref{Cinv} implies  that $K_{\mathcal F} \cdot C_1=-1$ and $K_{\mathcal F} \cdot C_i=0$ for $i=2,\dots,m$. Let $\varepsilon'\colon X\to X'$ be the birational morphism which contracts $C_2,\dots,C_m$. Then there exists a Cartier divisor $K$ on $X'$ such that $K_{\mathcal F}=\varepsilon'^{\ast} K$. We have $K\cdot \varepsilon'_*C_1=K_{\mathcal F}\cdot C_1=-1$ and in particular $\varepsilon'_*C_1$ is contained in the support of the negative part of the Zariski decomposition of $K$. Thus, $C_1,\dots,C_m$ are also contained in the  support of the negative part of the Zariski decomposition of $K_{\mathcal F}$ and the claim follows. 

We now assume that $\mathcal F$ is not relatively minimal.
Let $\varepsilon\colon Y\rightarrow (G\times F)/\Gamma$
be the minimal resolution.
By Remark \ref{impl2},
there exists a birational morphism $\nu\colon X\rightarrow Y$. By (1), it follows that $\nu$ is not an isomorphism. 
We prove the claim by induction on the number of blow-ups in $\nu$.
We can factor $\nu=\mu\circ \nu_1$ where $\nu_1\colon X\rightarrow X_1$ is the contraction of a single $\mathcal F$-exceptional curve $E_1$. 
Let $\mathcal{F}_1$ be the induced foliation on $X_1$ and 
let $K_{\mathcal{F}_1}=P_1+N_1$ be the Zariski decomposition of $K_{\mathcal F_1}$.
By induction, the support of $N_1$ is the exceptional locus of $\varepsilon\circ \mu$.
Since $\mathcal{F}_1$ has reduced singularities, there exists $\alpha\geq 0$ such that
$K_{\mathcal F}=\nu_1^{\ast}K_{\mathcal{F}_1}+\alpha E_1$. Thus, it is enough to show that $E_1$ is contained in the support of the negative part of the Zariski decomposition of $K_{\mathcal F}$.
If $\alpha >0$, then the claim follows immediately. If $\alpha=0$, then $K_{\mathcal F}=\nu_1^{\ast}K_{\mathcal{F}_1}$
and the support of the negative part of the Zariski decomposition of $K_{\mathcal F}$ is the preimage of the support of $N_1$.
Since $K_{\mathcal F}\cdot E_1=0$, the centre of $E_1$ in $X_1$ is a singular point of $\mathcal F_1$ and in particular, it is contained in the intersection of two components of a fibre of the induced morphism $X_1\to D$. Since the fibres of $(G\times F)/\Gamma\to D=G/\Gamma$ are irreducible,  at least one of these components is $(\varepsilon\circ\mu)$-exceptional. Thus, the claim follows.

Finally, $(3)$ follows from the fact that any fibre of the morphism $(G\times F)/\Gamma\to D=F/\Gamma$ is irreducible and its strict transform on $X$ is the only  component in the fibre $\psi^*p$ which is not $\mathcal F$-invariant. 
\end{proof}

\medskip
\subsubsection{Foliations of Kodaira dimension one}
We conclude this section with the following characterisation of foliations with reduced singularities of Kodaira dimension one, due to  Mendes and McQuillan \cite{mendes00, mcq08}:
\begin{teo}\label{t_kod1}
Let $\mathcal F$ be a foliation with reduced singularities on a smooth surface and such that $\kappa(\mathcal F)=1$. Then $\mathcal F$ is one of the following: 
\begin{enumerate}
\item a Riccati foliation;
\item a turbulent foliation;
\item a foliation induced by a non-isotrivial elliptic fibration; or 
\item a foliation induced by an isotrivial fibration of genus $\geq 2$.
\end{enumerate}
\end{teo}

%
%
%

\section{Proof of the main results}
Let $(X_t,\mathcal{F}_t)_{t\in\Delta}$ be a family of foliations on surfaces with reduced singularities. Our goal is to prove that for any sufficiently large positive integer $m$, the plurigenera $h^0(X_t,mK_{\mathcal F_t})$ does not depend on $t\in \Delta$.
By Theorem \ref{bruinv}, we can analyse separately the three cases
\begin{enumerate}
\item[(a)] $\kappa(\mathcal{F}_t)=0$ for all $t\in \Delta$; 
\item[(b)] $\kappa(\mathcal{F}_t)=1$ for all $t\in \Delta$; and 
\item[(c)] $\kappa(\mathcal{F}_t)=2$ for all $t\in \Delta$.
\end{enumerate}

\subsection{Kodaira dimension 0} We begin with case $(a)$ above:
\begin{prop}\label{p_kod0}
Let $(X_t,\mathcal{F}_t)_{t\in\Delta}$ be a family of foliations such that $\kappa(\mathcal{F}_t)=0$ for any $t\in \Delta$.
Then for any positive integer $m$, the dimension $h^0(X_t,mK_{\mathcal{F}_t})$ does not depend on $t\in\Delta$.
\end{prop}
\begin{proof}
Let $K_{\mathcal{F}_t}=P_t+N_t$ be the Zariski decomposition of
$K_{\mathcal{F}_t}$. 
By Lemma \ref{lcm}, the Cartier index $m_0$ of $N_t$ does not depend on $t\in \Delta$.
By  Lemma \ref{l_kod0}, the Cartier divisor $m_0P_t$ is a torsion divisor. In particular, 
the torsion index l does not depend on $t\in \Delta$.
Hence, for all $t\in\Delta$ we have 
$$h^0(X_t, mK_{\mathcal{F}_t})=\left\{
\begin{array}{ll}
0& \;\;\;{\rm if}\;\;\; m_0 l\;\;\;{\rm does}\;\;\;{\rm not}\;\;\;{\rm divide}\;\;\; m;\\
1& \;\;\;{\rm otherwise}.
\end{array}
\right.
$$
Thus, the claim follows.
\end{proof}

\subsection{Kodaira dimension 1} We now consider case $(b)$ above.  Given a family of foliations $(X_t,\mathcal{F}_t)_{t\in\Delta}$  of Kodaira dimension $1$, we first prove that the foliations $\mathcal F_t$
 are all of the same type:
either $\mathcal{F}_t$ is Riccati for all $t$, or turbulent for all $t$,
or induced by a non-isotrivial elliptic fibration  for all $t$, or by an isotrivial fibration of curves of genus greater than 1 for all $t$ (see Proposition \ref{p_kod1}).
Then,  we show that,  
for any $t\in \Delta$, there exists a fibration $\varphi_t\colon X_t\rightarrow C_t$ onto a curve $C_t$ 
such that $\varphi_t^{\ast}\delta_t$ is the positive part of the Zariski decomposition of $K_{\mathcal F_t}$ for some 
 ample $\mathbb{Q}$-divisor $\delta_t$ on $C_t$ such that  the degree of $\delta_t$ and the genus of $C_t$ do not depend on $t\in \Delta$.
Finally, for every such case, we show the invariance of the plurigenera  $h^0(X_t,mK_{\mathcal F_t})$ for a sufficiently large positive integer $m$.

We begin with the following basic and more general result:
\begin{lem}\label{l_kod1}
Let $\pi\colon\mathcal{X}\rightarrow\Delta$ be a smooth family of surfaces.  For any $t\in \Delta$, let $X_t=\pi^{-1}(t)$.
Let $P$ be a nef $\mathbb Q$-divisor such that 
$P\vert_{X_t}$ is of Kodaira dimension $1$ for any $t\in \Delta$ and let $P_t=P\vert_{X_t}$.

Then, for any $t\in \Delta$, there exists a fibration $\varphi_t\colon X_t\to C_t$ onto a smooth curve $C_t$ and an ample $\mathbb Q$-divisor $\delta_t$ on $C_t$ such that $P_t=\varphi_t^{\ast}\delta_t$. 

Further, if $P_t\cdot K_{X_t}\le 0$ for some $t\in \Delta$ then   the degree of $\delta_t$, the genus of the general fibre of $\varphi_t$ and the genus of $C_t$ do  not depend on $t\in \Delta$. 
\end{lem}
\begin{proof}
 By assumption, we have that $P_t$ is semi-ample for all $t\in \Delta$ and there exists a fibration $\varphi_t\colon X_t\to C_t$ onto a smooth curve $C_t$ and an ample $\mathbb Q$-divisor $\delta_t$ on $C_t$ such that $P_t=\varphi_t^{\ast}\delta_t$. 


Pick $s\in \Delta$ and let  
$F$ be a general smooth fibre of $\varphi_s$. 
We first assume  that $K_{X_s}\cdot P_s<0$. Then, for all $t\in \Delta$, the general fibre of $\varphi_t$ is a smooth rational curve and, in particular, $K_{X_s}\cdot F=-2$. Thus, for all $t\in \Delta$, we have that  $P\cdot K_{\mathcal X}\cdot X_t=-2\deg \delta_t$, which implies that the degree of $\delta_t$ does not depend on $t\in \Delta$.  Further, for all $t\in\Delta$, we have $R^1\varphi_{t\ast}\mathcal O_{X_t}=0$  and by the Leray spectral sequence we have that 
$$g(C_t)=h^1(C_t,\varphi_{t\ast}\mathcal{O}_{X_t})=h^1(X_t,\mathcal{O}_{X_t}).$$ 
Thus,  the claim follows. 

We now assume that $K_{X_s}\cdot P_s=0$. Then, the general fibre of $\varphi_t$ is an elliptic curve for all $t\in \Delta$.  Let us assume that there exists a $(-1)$-curve $E_s$ on $X_s$ which is contained in a fibre of $\varphi_s$. 
If $N_{E_s/X_s}$ denotes the normal bundle of $E_s$ in $X_s$ then 
$h^1({E_s},N_{{E_s}/X_{s}})=0$
and  by \cite[Theorem 1]{kodaira63}, after possibly shrinking $\Delta$, there exists a smooth surface $E$ in $\mathcal X$ 
which intersects the fibres of $\pi$ transversally and  such that $E|_{X_s}=E_s$. Thus, 
there exists a birational morphism $\varepsilon \colon \mathcal X\to \mathcal X'$ which contracts $E$ and a smooth morphism $\pi'\colon \mathcal X'\to \Delta$ such that $\pi=\pi'\circ\varepsilon$ (see also \cite[Chapt. 1, Prop. 1.20]{fm94} for a similar argument).
Let $Q=\varepsilon_*P$. Since $P_s\cdot E_s=0$, it follows that $P=\varepsilon^*Q$ and therefore we may replace $\mathcal X$ by $\mathcal X'$ and $P$ by $Q$. Thus, after finitely many steps, we may assume that there are no $(-1)$-curves on $X_s$ which are contained in a fibre of $\varphi_s$. Since $X_s$ is relatively minimal, it follows that $K^2_{X_s}=0$. Thus, $K^2_{X_t}=0$ for all $t\in \Delta$ and $X_t$ is also relatively minimal.

Let $B_{C_t}$ be the discriminant of $\varphi_t$ (cf. \eqref{e_disc}), let $M_{C_t}$ be the moduli part of $\varphi_t$ and let  $B'_{C_t}$ as in \eqref{e_disc2}, for all $t\in \Delta$. As in \cite[Chapt. I, Prop. 7.1]{fm94}, if $\varphi_s\colon X_s\to C_s$ has $k$ multiple fibres, of multiplicities $m_1,\dots,m_k$ for some $s\in \Delta$, then every surface $X_t$ has exactly $k$ multiple fibres of multiplicities $m_1,\dots,m_k$. 
By \cite[Chapt.  V, Prop.12.2]{bhpv04}, we have that $\deg (M_{C_t}+B'_{C_t})=\chi(X_t,\mathcal O_ {X_t})$ does not depend on $t\in \Delta$. Thus, it follows that the degree of  $M_{C_t}+B_{C_t}$ does not depend on $t\in \Delta$. 
Since $K_{X_t}=\varphi^*_t(K_{C_t}+M_{C_t}+B_{C_t})$,  invariance of plurigenera implies that the genus of $C_t$ is constant.


We now show that, after possibly shrinking $\Delta$, the divisor $P$ is semi-ample over $\Delta$, i.e. there exists a family of curves $\rho\colon \mathcal C\to \Delta$ 
and a morphism $q\colon \mathcal X\to  \mathcal C$ such that $\pi$ factors through $q$ and $ P\equiv q^*A$ for some 
 ample $\mathbb Q$-divisor $A$ on $\mathcal C$. 
To this end, we distinguish three  cases. If  $\kappa(X_t)>0$ then \cite[Chapt. I, Thm. 7.11 (iii)] {fm94} implies the existence of the family $\rho\colon \mathcal C\to \Delta
$ and the morphism $q\colon\mathcal X\to \mathcal C$ as above. Let $A$ be any relatively ample $\mathbb Q$-divisor on $\mathcal C$ and let $P'=q^*A$. Then, after possibly rescaling $A$, we may assume that $P|_s \equiv P'|_s$. Thus, $(P-P')\cdot H\cdot X_t=0$ for any ample divisor $H$ on $\mathcal X_t$ and for all $t\in \Delta$. It follows that, after possibly shrinking $\Delta$, we have that $P\equiv P'$ and the claim follows. 

Thus, we may assume that $\kappa(X_t)\le 0$ and, in particular, $g(C_t)\le 1$ for all $t\in \Delta$.
We first assume that $g(C_t)=0$ for all $t\in \Delta$. 
After possibly shrinking $\Delta$, we may assume that the stable base locus of $P$ is contained in the fibre $X_0=\pi^{-1}(0)$ and in particular, if $\mathcal X^{\ast}=\mathcal X\setminus X_0$, then $P|_{\mathcal X^{\ast}}$ is semi-ample. 
 Thus, there exists a factorization 
$$\mathcal{X^{\ast}}\xrightarrow{\varphi}\mathcal{C^{\ast}}\xrightarrow{\nu}\Delta\setminus \{0\}$$
where, for each $t\in \Delta\setminus \{0\}$, we have $C_t=\nu^{-1}(t)$. 
We may  assume that $P\ge 0$. Let  $x\in X_0$  be a general point. After possibly shrinking $\Delta$, we may assume that there exists a one-dimensinoal subvariety $\Gamma\subseteq \mathcal X$ such that the induced morphism $\pi|_{\Gamma}\colon \Gamma\to \Delta$ is an isomorphism and $\Gamma$ meets $X_t$ transversally for each $t\in \Delta$. Let $\Gamma^*=\Gamma\setminus X_0$ and let $R^*=\varphi^{-1}(\varphi(\Gamma^*))$. Then, for each $t\in \Delta\setminus\{0\}$, we have that $R^*\cap X_t$ is a fibre of $\varphi_t$. Let $R$ be the closure of $R^*$ in $\mathcal X$. After possibly shrinking $\Delta$, we may assume that $R$ does not intersect the support of $P$. 
There exists a positive rational number $\alpha$ such that if $Q=\alpha R$, we have that $P|_{X_0}\equiv Q|_{X_0}$ for some $\mathbb Q$-divisor $Q\ge 0$ whose support does not intersect the support of $P$. As above, it follows that after possibly shrinking $\Delta$ further, we have that $P\equiv Q$. 
Moreover, since $g(C_t)=0$,  there exists a positive integer $m$ such that  $mP|_{X_t}\sim m Q|_{X_t}$ for all $t\in \Delta$. In particular,  $mP\sim mQ$. Since the support of $P$ and the support of $Q$ do not intersect, it follows that the relative Iitaka fibration associated to $P$ is a proper morphism over $\Delta$, i.e.  $P$ is semi-ample over $\Delta$.  Thus,  the claim follows. 

We now assume that $g(C_t)=1$, and in particular $\kappa(X_t)=0$ for all $t\in \Delta$.  By the Enrique-Kodaira classification of algebraic surfaces, it follows that $X_t$ is either hyperelliptic for all $t\in \Delta$ or an abelian surface for all $t\in \Delta$. In the first case, \cite[pag. 130, Remark 2)] {fm94} implies the existence of the family $\rho\colon \mathcal C\to \Delta$ and we can conclude as above. Let us assume now that $X_t$ is an abelian surface for all $t\in \Delta$.   
By \cite[Thm 6.14]{mumford94} (see also \cite[pag. 89]{FC90}), we may assume that $\mathcal X$ is an abelian scheme over $\Delta$. In particular, there exists a morphism  $\beta \colon \mathcal X\to \mathcal X$ over $\Delta$, such that its restriction to $X_t$ coincides with the inverse morphism on $X_t$ for all $t\in \Delta$ and, for each section $\gamma\colon \Delta\to \mathcal X$ of $\pi$, we may   define a morphism $t_\gamma\colon \mathcal X\to \mathcal X$ over $\Delta$, whose restriction on $X_t$ coincides with the translation by $\gamma(t)$ on $X_t$ for all $t\in \Delta$.  
Thus, the theorem of squares (e.g. \cite[Theorem 2.3.3]{BL04}) implies that if $Q$ is any Cartier divisor on $\mathcal X$ and if $\gamma'=\beta\circ\gamma$, then 
$2Q\sim t^*_\gamma Q\otimes t^*_{\gamma'} Q$. In particular, if $Q$ is effective, then after possibly shrinking $\Delta$, we have that $Q$ is semi-ample over $\Delta$.   Thus, also in this case, the claim follows. 

We now show that  the degree of $\delta_t$ does not depend on $t\in \Delta$. We have  that 
 $ P\equiv q^*A$ for some 
relatively ample $\mathbb Q$-divisor $A$ on $\mathcal C$. 
Note that since the fibres of $\pi$ are reduced,  also the fibres of $\rho$ are reduced. Thus, we have that  
$$\deg \delta_t=\deg A|_{\rho^{-1}(t)}$$
does not depend on $t\in \Delta$, as claimed.  
\end{proof}

We now show that the type of the foliation is preserved. A similar argument also appeared  in \cite[pag. 130]{brunella01}.
\begin{prop}\label{p_kod1}
Let $(X_t,\mathcal{F}_t)_{t\in\Delta}$ be a family of foliations such that $\kappa(\mathcal{F}_t)=1$
for any $t\in \Delta$.
Then all the $\mathcal{F}_t$ are of the same type, i.e. there are $4$ possibilities:
\begin{enumerate}
\item $\mathcal{F}_t$ is a Riccati  foliation transverse to a fibration $\varphi_t\colon X_t\to C_t$  for all $t\in \Delta$ (cf. Definition \ref{d_riccati});
\item $\mathcal{F}_t$ is a turbulent  foliation transverse to a fibration $\varphi_t\colon X_t\to C_t$  for all $t\in \Delta$ (cf. Definition \ref{d_turbulent});
\item $\mathcal{F}_t$ is induced by a non-isotrivial elliptic fibration  $\varphi_t\colon X_t\to C_t$ for all $t\in \Delta$;
\item $\mathcal{F}_t$ is not Riccati or turbulent and it is induced by an isotrivial fibration of genus $\geq 2$ for all $t\in \Delta$.
\end{enumerate}
Further, in the cases $(1)$, $(2)$ and $(3)$  the  genus $g(C_t)$ of the curve $C_t$ does not depend on $t\in \Delta$. 
\end{prop}

\begin{proof}
We first prove that the class of the foliation  $\mathcal F_t$ does not depend on $t\in \Delta$. 
Fix $s\in\Delta$. By Lemma \ref{Fecc}, we may assume that $\mathcal F_s$ is relatively minimal. Let $N\ge 0$ be the $\mathbb Q$-divisor on $\mathcal X$ whose existence is guaranteed by Proposition \ref{zardecfol}, such that, if  we denote $P=K_{\mathcal F}-N$, $P_t=P|_{X_t}$ and $N_{t}=N|_{X_t}$ then $K_{\mathcal F_t}=P_t+N_t$ is the Zariski decomposition of $K_{\mathcal F_t}$ for all $t\in \Delta$. By Theorem \ref{structneg}, we have that $\lfloor N_s\rfloor=0$. 
 By Proposition \ref{zardecfol}, after possibly shrinking $\Delta$, we may assume that
  each irreducible component of $N$ meet the surfaces $X_t$ transversally in a rational curve and that $\lfloor N_t\rfloor=0$ for all $t\in \Delta$.

Let $t\in \Delta$. 
By Theorem \ref{t_kod1}, $\mathcal F_t$ is one of the following: 
\begin{enumerate}
\item a Riccati foliation;
\item a turbulent foliation;
\item a foliation induced by a non-isotrivial elliptic fibration; or 
\item a foliation which is not Riccati or turbulent and it is induced by an isotrivial fibration of genus $\geq 2$.
\end{enumerate}

For all $t\in \Delta$, let  $\varphi_t\colon X_t\to C_t$ be the Iitaka fibration of $P_t$.
Note that, by \eqref{cantrasv}, if $\mathcal F_t$ is of type $(1)$, $(2)$  then it is transverse to $\varphi_t$ and 
by (2) of Proposition \ref{cbf},  if $\mathcal F_t$ is of type $(3)$ then it is induced by $\varphi_t$.
If $\mathcal F_t$ is of type $(4)$ and $F$ is the general fibre of $\varphi_t$,
since $\varphi_t$ is the Iitaka fibration of $K_{\mathcal F_t}$, it follows that  $K_{\mathcal F_t}\cdot F=0$. 
Since any $\mathcal F_t$-invariant curve passing through the general point of $X_t$ is of genus greater than one, Proposition \ref{Cinv} implies that  $F$ is not $\mathcal F_t$-invariant and therefore Proposition \ref{Cnotinv} implies that $\mathcal F_t$ is transverse to $\mathcal \varphi_t$.  By assumption $\mathcal F_t$ is not Riccati or turbulent, and in particular the genus of $F$ is greater than one. 

The intersection $K_{X_t}\cdot P\vert_{X_t}=K_{\mathcal X}\cdot P\cdot X_t$ does not depend on $t\in \Delta$.
If $K_{X_t}\cdot P\vert_{X_t}<0$ for all $t\in \Delta$ then $\mathcal F_t$ is of type $(1)$ and 
if $K_{X_t}\cdot P\vert_{X_t}>0$ for all $t\in \Delta$ then $\mathcal F_t$ is of type $(4)$.
Thus, we may assume that $\varphi_t$ is an elliptic fibration and $\mathcal F_t$ is either turbulent or induced by a non-isotrivial elliptic fibration. 

Let $B_{C_t}$ be the discriminant of $\varphi_t$ (cf. \eqref{e_disc}), let $M_{C_t}$ be the moduli part of $\varphi_t$ and let  $B'_{C_t}$ as in \eqref{e_disc2}, for all $t\in \Delta$. 
By \cite[Chapt.  V, Prop.12.2]{bhpv04}, we have that $\deg (M_{C_t}+B'_{C_t})=\chi(X_t,\mathcal O_ {X_t})$ does not depend on $t\in \Delta$.

For any $t\in \Delta$ and for any $p\in {\rm Supp} B'_{C_t}$, (1) and (3) of Corollary  \ref{c_turb} imply that there exists a prime divisor $D\subseteq \varphi_t^{-1}(p)\cap {\rm Supp } N$ which does not compute the log canonical threshold at $p$ (cf. Remark  \ref{r_lct}). 
We claim that the number of points in ${\rm Supp} B'_{C_t}$ does not depend on $t\in \Delta$.
Indeed, if $\mathcal{C},\mathcal{C}'$ are connected components of ${\rm Supp}N$
such that $\mathcal{C}\vert_{X_t},\mathcal{C}'\vert_{X_t}$ are contained in two different fibres of $\varphi_t$
for some $t\in \Delta$, then it is sufficient to show that 
$\mathcal{C}\vert_{X_t},\mathcal{C}'\vert_{X_t}$ are chains of rational curves which are contained in two different fibres of $\varphi_t$  for all $t\in \Delta$.
Indeed, assume by contradiction that there exists $t_0\in \Delta$
and $\mathcal{C}_1,\dots, \mathcal{C}_h$,
$\mathcal{C}_{h+1},\dots, \mathcal{C}_k$
connected components of ${\rm Supp}N$
such that 
\begin{itemize}
\item $\mathcal{C}_i\vert_{X_t}$ is contained in $F_1^t$
for $i=1,\dots, h$ and
\item $\mathcal{C}_i\vert_{X_t}$ is contained in $F_2^t$ for $i=h+1,\dots, k$
\end{itemize}
where $F_1^t,F_2^t$ are fibres of $\varphi_t$ and $F_1^t\neq F_2^t$
for $t\neq t_0$ and $F_1^{t_0}= F_2^{t_0}$.
Since $\lfloor N_t\rfloor=0$ for all $t\in \Delta$, Lemma \ref{D_tan} implies that the chains
$\mathcal{C}_i\vert_{X_t}$ contain all the reduced components of $F_1^t$ for $i=1,\dots, h$ for $t\neq t_0$. It follows that the same condition is preserved for $t=t_0$. 
By (1) of Corollary \ref{c_turb}, the fibre $F_2^t$ is not a multiple fibre for all $t\neq t_0$ and 
since $\varphi_t$ is a fibration by elliptic curves,
 there exists a reduced component of $F_2^t$, for all $t\neq t_0$.
Thus, there exists a component of a chain $\mathcal{C}_i$
with $i=h+1,\dots, k$ that degenerates to a rational curve with multiple coefficient,
which is a contradiction because, by Proposition \ref{zardecfol}, the chains $\mathcal{C}_i$ meet the fibres of $\pi$
transversally. Thus, it follows that the number of the points in the support of $B'_{C_t}$ does not depend on $t\in \Delta$. 

 Moreover, by the classification of singular fibres of an elliptic fibration, there exists a reduced component $E$ of $\varphi_t^*p$ which is not contained in the exceptional locus of the relative minimal fibration $X_t\to X'_t$ associated to $\varphi_t\colon X_t\to C_t$. Thus, (2) of  Corollary \ref{c_turb} implies that $\deg B'_{C_t}$ does not depend on $t\in \Delta$.  In particular, $\deg M_{C_t}$ does not depend on $t\in \Delta$.
 Remark \ref{r_turbulent} implies that $\deg M_{C_t}=0$ if and only if $\mathcal F_t$ is turbulent. 
Thus, either $\mathcal F_t$ is turbulent for all $t\in \Delta$ or $\mathcal F_t$ is induced by a non-isotrivial elliptic fibration for all $t\in \Delta$. 

Finally, Lemma \ref{l_kod1} implies that  in the cases $(1)$, $(2)$ and $(3)$, the genus of the curve $C_t$ does not depend on $t\in \Delta$. 
%
%
%
%
\end{proof}

\begin{prop}\label{p_kod1a}
Let $(X_t,\mathcal{F}_t)_{t\in\Delta}$ be a family of foliations induced by non-isotrivial elliptic fibrations $\varphi_t\colon X_t\to C_t$ over a curve $C_t$ of genus $g$.
Then for any sufficiently large positive intger $m$, the dimension $h^0(X_t,mK_{\mathcal{F}_t})$
does not depend on $t\in \Delta$.
\end{prop}

In Example \ref{e_elliptic}, we show that in general the claim does not hold if $m=1$. 

\begin{proof}
We want to show  that for any sufficiently large positive integer $m$, the dimension $h^0(X_t,mK_{\mathcal{F}_t})$
is locally constant.
Fix $s\in\Delta$. By Lemma \ref{Fecc}, we may assume that $\mathcal F_s$ is relatively minimal.
Let $N\ge 0$ be the $\mathbb Q$-divisor on $\mathcal X$ whose existence is guaranteed by Proposition \ref{zardecfol}, such that, if  we denote $P=K_{\mathcal F}-N$, $P_t=P|_{X_t}$ and $N_{t}=N|_{X_t}$ then $K_{\mathcal F_t}=P_t+N_t$ is the Zariski decomposition of $K_{\mathcal F_t}$ for all $t\in \Delta$. By Theorem \ref{structneg}, we have that $\lfloor N_s\rfloor=0$. 
 By Proposition \ref{zardecfol}, after possibly shrinking $\Delta$, we may assume that
  each irreducible component of $N$ meet the surfaces $X_t$ transversally in a rational curve and that $\lfloor N_t\rfloor=0$ for all $t\in \Delta$.


For any $t\in \Delta$, let $\varepsilon_t\colon X_t \to X'_t$ be the minimal elliptic fibration associated to $\varphi_t$ and let $\varphi'_t\colon X'_t\to C_t$ be the induced morphism. 
Let $B_{C_t}$ and $M_{C_t}$ be the discriminant and the moduli part in the canonical bundle formula of $\varphi_t$ (cf. \eqref{e_disc}). Then,   Lemma \ref{cbf} implies that 
$$K_{X_t}=\varphi_t^{\ast}(K_{C_t}+M_{C_t}+B_{C_t})+E_t\qquad\text{and}\qquad P_t=\varphi_t^{\ast}(M_{C_t}),$$
where $E_t\ge 0$ is $\varepsilon_t$-exceptional for each $t\in \Delta$. 

Then Lemma \ref{l_kod1} implies that 
$$\deg(K_{C_t}+M_{C_t}+B_{C_t})\qquad\text{and}\qquad\deg(M_{C_t})$$ do not depend on $t\in \Delta$.
By Proposition \ref{p_kod1} the genus of $C_t$ does not depend on $t\in \Delta$. Thus,  also $\deg B_{C_t}$ does not depend on $t\in \Delta$.

For any $t\in \Delta$, let $B'_{C_t}$ be the $\mathbb Q$-divisor on $C_t$ defined as in \eqref{e_disc2}. By \eqref{fracpart}, we have  $${\rm Supp}\{M_{C_t}\}={\rm Supp}B'_{C_t}.$$

As in the proof of Proposition \ref{p_kod1}, there exists $c_1,\dots,c_k\in (0,1)$ and, 
for any $t\in \Delta$, there exist distinct $p_1^t,\dots,p_k^t\in C_t$
such that $$B'_{C_t}=\sum_{i=1}^k (1-c_i)p_i^t.$$
Since $M_{C_t}+B'_{C_t}$ is integral,
 $$\{M_{C_t}\}=\sum_{i=1}^k c_i p_i^t.$$
Thus, for any positive integer $m$ the degree of $\lfloor mM_{C_t}\rfloor$, does not depend on $t\in \Delta$ and we may choose $m$  sufficiently large so that  
$\deg \lfloor mM_{C_t}\rfloor\ge 2g-1$. 
For all $t\in \Delta$, we have 
$$h^0(X_t, \mathcal O_{X_t}(\lfloor mP_t\rfloor))=h^0(C_t, \mathcal O_{C_t}(\lfloor mM_{C_t}\rfloor))=\deg \lfloor mM_{C_t}\rfloor+1-g(C_t)$$
and the claim follows. 
\end{proof}

\begin{prop}\label{p_kod1b}
Let $(X_t,\mathcal{F}_t)_{t\in\Delta}$ be a family of Riccati or turbulent foliations  transverse to the fibration $\varphi_t\colon X_t\to C_t$ for all $t\in \Delta$.
Then for any positive integer $m$, the dimension $h^0(X_t,mK_{\mathcal{F}_t})$
does not depend on $t\in \Delta$.
\end{prop}
\begin{proof}
We want to show that $h^0(X_t,mK_{\mathcal{F}_t})$ is locally constant. 
Fix $s\in\Delta$. By Lemma \ref{Fecc}, we may assume that $\mathcal F_s$ is relatively minimal. Let $N\ge 0$ be the $\mathbb Q$-divisor on $\mathcal X$ whose existence is guaranteed by Proposition \ref{zardecfol}, such that, if  we denote $P=K_{\mathcal F}-N$, $P_t=P|_{X_t}$ and $N_{t}=N|_{X_t}$ then $K_{\mathcal F_t}=P_t+N_t$ is the Zariski decomposition of $K_{\mathcal F_t}$ for all $t\in \Delta$. By Theorem \ref{structneg}, we have that $\lfloor N_s\rfloor=0$. 
 By Proposition \ref{zardecfol}, after possibly shrinking $\Delta$, we may assume that
$N$ meets the fibres of $\pi$ transversally in chains of rational curves and that $\lfloor N_t\rfloor=0$ for all $t\in \Delta$.
By \eqref{cantrasv}, there exist effective divisors $D^t_{\tan}$ and $D_{\varphi_t}$ contained in the fibres of $\varphi_t$ which are $\mathcal F_t$-invariant and such that  $$K_{\mathcal{F}_t}=\varphi_t^{\ast}K_{C_t}+D_{\tan}^t+D_{\varphi_t}.$$
As in \eqref{Dtan}, there exists a $\mathbb{Q}$-divisor $\theta_t$ on $C_t$
such that $D^t_{\tan}+D_{\varphi_t}=\varphi_t^{\ast}\theta_t+N_t$ and $P_t=\varphi_t^*(K_{C_t}+\theta_t)$, for any $t\in \Delta$ .
By Lemma \ref{l_kod1},
the degree of $K_{C_t}+\theta_t$ is constant and, by Proposition \ref{p_kod1}, the genus of  the curve $C_t$ does not depend on $t\in \Delta$. It follows that also the degree of $\theta_t$ is constant.

For all $t\in \Delta$, we may write 
$$\theta_t=\sum_{p\in C_t} \theta_{p,t}p.$$

We may assume $s=0$. 
%
We want to show that the components of $\varphi_t^{\ast}p$
such that the coefficient $\theta_{p,t}$ is not integral do not meet as $t$ approaches $0$.
Let  $p\in C_{t}$ be such that $\theta_{p,t}$ is rational but not integral, for some $t\in \Delta$. 
We claim that either $\varphi_{t}^{\ast}p$ is a multiple fibre whose support is a smooth curve  $\tilde{F}$ of genus one
or it contains a reduced component.
Indeed if $\mathcal F_t$ is a Riccati foliation, then the claim follows from the fact that every fibre of $\varphi_t$ admits a reduced component. On the other hand, if $\mathcal F_t$ is a turbulent foliation, then the claim follows from (4) of Corollary \ref{c_turb} and the fact that, by classification of singular fibres of elliptic fibrations,  any singular fibre, which is not a multiple fibre, contains a reduced component.

 As in \cite[Chapt. I, Prop. 7.1]{fm94}, if $\varphi_0\colon X_0\to C_0$ has $k$ multiple fibres, of multiplicities $m_1,\dots,m_k$, then every surface $X_t$ has exactly $k$ multiple fibres of multiplicities $m_1,\dots,m_k$. 
It follows that multiple fibres cannot meet.
We now assume that $\varphi_t^{\ast}p$ is not a multiple fibre and it contains  a reduced component.
By Lemma \ref{D_tan},
${\rm Supp}N_t$ contains all the reduced components of $\varphi_t^{\ast}p$. Thus, 
since  $N$ meets $X_t$ transversally  for all $t\in \Delta$, it follows that 
  the components of $\varphi_t^*p$ such that the coefficient $\theta_{p,t}$ is not integral do not meet as $t$ approaches  0. 
  Indeed the same argument as in Proposition \ref{p_kod1} goes through 
  because since $\varphi_t$ is a fibration by rational or elliptic curves,
 every fibre has a reduced component.

By Proposition \ref{p_kod1}, the genus $g(C_t)$ does not depend on $t\in \Delta$. Thus,
$$h^0(X_t,mK_{\mathcal{F}_t})=h^0(X_t,\lfloor mP_t\rfloor)=h^0(C_t, mK_{C_t}+\lfloor m\theta_t\rfloor)=
\deg (mK_{C_t}+\lfloor m\theta_t\rfloor) +1-g(C_t)$$
does not depend on $t\in \Delta$ and the claim follows.
\end{proof}

\begin{prop}\label{p_kod1d}
Let $(X_t,\mathcal{F}_t)_{t\in\Delta}$ be a family of foliations which are not Riccati or turbulent and they are induced by isotrivial fibrations
of genus $g\geq2$.
Then for any positive integer $m$, the dimension $h^0(X_t,mK_{\mathcal{F}_t})$
does not depend on $t\in \Delta$.
\end{prop}
\begin{proof}
We want to show that $h^0(X_t,mK_{\mathcal{F}_t})$ is locally constant. 
Fix $s\in\Delta$. 
By Lemma \ref{Fecc}, we may assume that $\mathcal F_0$ is relatively minimal.
Let $N\ge 0$ be the $\mathbb Q$-divisor on $\mathcal X$ whose existence is guaranteed by Proposition \ref{zardecfol}, such that, if  we denote $P=K_{\mathcal F}-N$, $P_t=P|_{X_t}$ and $N_{t}=N|_{X_t}$ then $K_{\mathcal F_t}=P_t+N_t$ is the Zariski decomposition of $K_{\mathcal F_t}$ for all $t\in \Delta$. By Theorem \ref{structneg}, we have that $\lfloor N_0\rfloor=0$. 
 By Proposition \ref{zardecfol}, after possibly shrinking $\Delta$, we may assume that
  each irreducible component of $N$ meet the surfaces $X_t$ transversally in a rational curve and that $\lfloor N_t\rfloor=0$ for all $t\in \Delta$.

By Remark \ref{impl2}, for all $t\in \Delta$, there exists a proper birational morphism  
$$\alpha_t\colon X_t\to W_t := (G_t\times F_t)/\Gamma_t$$
 where $F_t$ and $G_t$ are smooth curves such that the genus of $F_t$ is greater than one and $\Gamma_t$ is a finite group acting on $F_t$ and $G_t$. By  Theorem \ref{serr41}, for any $t\in \Delta$, 
there exist two fibrations $\varphi_t\colon X_t\to C_t$ and $\psi_t\colon X_t\to D_t$ over the curves $C_t=F_t/\Gamma_t$ and $D_t=G_t/\Gamma_t$ with fibres having support with only normal crossing singularities and such that $\mathcal F_t$ is induced by $\psi_t$ and it is transverse to $\varphi_t$. Since $\mathcal F_t$ is not Riccati nor turbulent, also the curve $G_t$ is of genus greater than one for all $t\in \Delta$. 

We now show that  $\mathcal{F}_t$ is relatively minimal  for all $t\in\Delta$.
Indeed,  by $(2)$ of Lemma \ref{overpoint}, the support of $N_t$ coincides with the exceptional locus of $\alpha_t$.
By $(1)$ of Lemma \ref{overpoint}, the  exceptional locus of $\alpha_t$ contains a $(-1)$-curve if and only if the foliation is not relatively minimal. Since $\mathcal F_s$ is relatively minimal, it follows that the support of $N_s$ is a union of Hirzebruch-Jung strings and, in particular, the support of $N_s$ does not contain any $(-1)$-curve. Therefore the support of $N_t$ does not contain any $(-1)$-curve for all $t\in \Delta$. Thus, $\mathcal{F}_t$ is relatively minimal  and $X_t$ is the minimal resolution of $W_t$ for all $t\in \Delta$.

 Let $E_1,\dots, E_k$ be the irreducible components of $N$ and, for any $t\in \Delta$, let 
$$E_i^t={E_i}|_{X_t}.$$ 
By Theorem \ref{structneg}, $E_i^t$ is a smooth rational curve. In particular, the self-intersection $(E_{i}^t)^2$ is a negative number which does not depend on $t\in \Delta$. 
Since, by (2) of Lemma \ref{overpoint}, the exceptional locus of  $\alpha_t$ coincides with the support of $N_t$, it follows that  there exists an effective $\mathbb Q$-divisor $\Theta$ on $\mathcal X$, whose coefficients are contained in $(0,1)\cap \mathbb Q$ and whose support is contained in the support of $N$, such that 
$$K_{X_t}=\alpha_t^*K_{W_t}-\Theta|_{X_t}$$
for any $t\in \Delta$. In particular, $(\mathcal X,\Theta)$ is Kawamata log terminal, and since $K_{W_t}$ is ample for all $t\in \Delta$, it follows by the base point free theorem \cite[Thm. 3.24]{KM98} (see also \cite[Example 2.17]{KM98})  that $K_{\mathcal X} + \Theta$ is semi-ample over $\Delta$, i.e.  if we define 
$$R(\mathcal X/\Delta,K_{\mathcal X}+\Theta)=\bigoplus _{m\ge 0} \pi_*\mathcal O_{\mathcal X}(\lfloor m(K_{\mathcal X}+\Theta)\rfloor )$$ and  
\begin{equation}\label{e_W}
\mathcal W={\rm Proj}~ R(\mathcal X/\Delta,K_{\mathcal X}+\Theta)
\end{equation}
 denotes the log canonical model of $K_{\mathcal X}+\Theta$ over $\Delta$,  then, there exists a proper birational morphism $\alpha\colon \mathcal X\to \mathcal W$ which defines a factorization of $\pi$ and such that, if we denote   $\Delta^{\ast}=\Delta\backslash\{0\}$, then, after possibly shrinking $\Delta$, the induced morphism $X_t\to \alpha(X_t)$ coincides with $\alpha_t$,  for all $t\in \Delta^*$. Note that, by the Negativity Lemma,  we have that $P=\alpha^*Q$, where $Q=\alpha_*P$ is a nef $\mathbb Q$-divisor on $\mathcal W$.




 We claim  that, after possibly shrinking $\Delta$, the $\mathbb Q$-divisor $P$ is numerically equivalent to a semi-ample over $\Delta$, i.e. there exists a family of curves $\mathcal C\to \Delta$ and a morphism $\varphi\colon \mathcal X\to \mathcal C$ which defines a factorization of $\pi$, and such that $P\equiv \varphi^*H$ for some relatively ample $\mathbb Q$-divisor $H$ on $\mathcal C$. To this end, we may and will perform a base change $\tau\colon \Delta\to \Delta$ which is totally ramified over the origin. Indeed, we may replace $\pi\colon \mathcal X\to \Delta$ by the family of surfaces $\pi'\colon\mathcal X'\to \Delta$ obtained by base change and $P$ by the pull-back of $P$ on $\mathcal X'$. We will do this, without changing notation. 

Lemma \ref{iitaka} and \eqref{e_kiso} imply that 
$\varphi_t$ is the Iitaka fibration of $K_{\mathcal F_t}$  for all $t\in \Delta$.
Let $Z=\cup E_i$ be the support of $N$. By $(2)$ of Lemma \ref{overpoint}, $Z|_{X_t}$ coincides with the exceptional divisor of $\alpha_t$.   Thus,  \eqref{e_kiso} and  Theorem \ref{serr41} imply that 
$$K_{X_t}-K_{\mathcal F_t} +Z|_{X_t}=K_{X_t} - K_{X_t/D_t}+D_{\psi_t} +Z|_{X_t} =  K_{X_t/C_t} -D_{\varphi_t}$$ 
and, therefore, Lemma \ref{iitaka} implies that 
$\psi_t$ is the Iitaka fibration of $(K_{\mathcal X} - K_{\mathcal F}+Z)\vert_{X_t}$ for all $t\in \Delta$.
It follows that, after possibly shrinking $\Delta$, we  may assume that, if $\mathcal X^{\ast}=\mathcal X\setminus X_0$,  
we have two factorisations
$$\mathcal{X}^{\ast}\xrightarrow{\varphi'}\mathcal{C}^{\ast}\rightarrow\Delta^{\ast}\qquad\text{and}\qquad\mathcal{X}^{\ast}\xrightarrow{\psi'}\mathcal{D}^{\ast}\rightarrow\Delta^{\ast}$$
such that for all $t\in \Delta^{\ast}$ the restrictions of $\varphi'$ and $\psi'$ to $X_t$ coincide with $\varphi_t$ and $\psi_t$.
After possibly shrinking $\Delta$ again, we may assume that there exist two families of curves $p'\colon F^*\to \Delta^*$ and $q'\colon G^*\to \Delta^*$ such that for all $t\in \Delta^*$, we have that $F_t$ is isomorphic to  $p'^{-1}(t)$ and $G_t$ is isomorphic to $q'^{-1}(t)$. Finally, we may assume that the group $\Gamma_t$ does not depend on $t\in \Delta^*$. We will denote it by $\Gamma$.    

Since the moduli functor of stable curves with an action of a finite group is 
 proper (e.g. see \cite{tuffery93}), after a base change $\tau\colon \Delta\rightarrow \Delta$ totally ramified over the origin, we may find two  families of curves $p\colon F\to \Delta$ and $q\colon G\to \Delta$ which extend the families $F^*$ and $G^*$ on the whole $\Delta$. Moreover, the group $\Gamma$  acts on $F$ and $G$ so that, for all $t\in \Delta^*$,  the action on $p^{-1}(t)$ and $q^{-1}(t)$ coincides with the action of $\Gamma_t$ on $F_t$ and $G_t$ respectively (see \cite[Theorem 3.1]{opstall06} for a similar argument). 
 
Let $\rho\colon \mathcal W'=(G\times_{\Delta}F)/\Gamma\to \Delta$. By abuse of notation, we continue to denote 
by $\mathcal W$, the pull-back of the  log canonical model over $\Delta$ defined in \eqref{e_W} after the base change performed above.   Note that $\mathcal W$ and $\mathcal W'$ are isomorphic over $\Delta^*$, and by the separateness property of the moduli functor of stable pairs (e.g. see \cite[Lemma 7.2]{hxu13}), it follows that $\mathcal W$ is isomorphic to $\mathcal W'$. In particular, there exists a proper birational morphism $\mathcal X\to \mathcal W'$, induced by $\alpha\colon \mathcal X\to \mathcal W$. 
Let $\mathcal C:=F/\Gamma$ and let $\nu\colon \mathcal C\to \Delta$ be the induced morphism. Then, there exists a proper morphism $\varphi\colon \mathcal X\to \mathcal C$ which factors through $\mathcal X\to \mathcal W'$ and, after possibly shrinking $\Delta$ again, we may assume that  there exists a relatively ample $\mathbb Q$-divisor $H$ over $\Delta$, and after rescaling $H$ we have $P\equiv \varphi^*H$. Thus, the claim follows. 

Since $C_0$ is the normalization of the curve $\nu^{-1}(0)\subseteq \mathcal C$,  we have that  $g(C_0)\le g(C_t)$ for all $t\in \Delta$ and, similarly,  we have that $g(D_0)\le g(D_t)$. Since \eqref{e_h1} implies that
$h^1(X_t,\mathcal{O}_{X_t})=g(C_t)+g(D_t)$, it follows that the genus of $C_t$ and the genus of $D_t$ do not depend on $t\in \Delta$.

Let $t\in \Delta$. By \eqref{cantrasv} there exists  effective divisors $D^t_{\tan}$ and $D_{\varphi_t}$ contained in fibres of $\varphi_t$ such that
 $$K_{\mathcal{F}_t}=\varphi_t^{\ast}K_{C_t} + D^t_{\tan}+D_{\varphi_t}.$$
 Thus, \eqref{e_kiso} and Theorem \ref{serr41} imply that 
  \begin{equation}\label{e_tan}
  Z\vert_{X_t}=D^t_{\tan}
  \end{equation}
where, as above,  $Z$ denotes the support of  $N$.

By \eqref{Dtan}, we have that $P_t=\varphi_t^{\ast}(K_{C_t}+\theta_t)$ for some $\mathbb Q$-divisor  $\theta_t\ge 0$ on $D_t$ such that $\deg (K_{C_t}+ \theta_t)$ does not depend on $t\in \Delta$ and $\theta_t$ is the largest $\mathbb Q$-effective divisor such that  $D^t_{\tan}+D_{\varphi_t}-\varphi_t^*\theta_t$ is effective. 
Since $g(C_t)$ is constant, we have that $\deg \theta_t = \deg \theta_0$ for all $t\in \Delta$.

We want to show that the points in the support of $\theta_t$  do not collide as $t$ approaches  $0$. 
For any $t\in \Delta$, the support of $\theta_t$ is contained in 
$\varphi_t({\rm Supp }(D^t_{\tan}+D_{\varphi_t}))$.
Multiple fibres with smooth support  do not collide
nor they meet components corresponding to singular fibres because the self-intersection of the reduced part of a multiple fibre is  nilpotent but non-trivial and the structure is preserved for all $t\in \Delta$.  Since the components of $N$ do not meet, (2) and (3) of Lemma \ref{overpoint} imply that the points in $\varphi_t({\rm Supp }N_t)$ do not collide as $t$ approaches  $0$. 
Thus,  \eqref{e_tan} implies the claim.

%
%
In particular,  there exist positive rational numbers $c_1,\dots,c_\ell$
such that for all $t\in \Delta$ there exist distinct points $p_{1,t},\dots,p_{\ell,t}$ which form irreducible curves on $\mathcal C$ and such that 
\begin{eqnarray}\label{gamma}
\theta_t&=&\sum c_i p_{i,t} \qquad \text{for all } t\in \Delta.
\end{eqnarray}
In particular, the degree $\deg\lfloor m\theta_t\rfloor$
does not depend on $t\in \Delta$. Thus,   
$$h^0(X_t,mK_{\mathcal{F}_t})=h^0(X_t,\lfloor mP_t\rfloor)=h^0(C_t, mK_{C_t}+\lfloor m\theta_t\rfloor)=
\deg (mK_{C_t}+\lfloor m\theta_t\rfloor) +1-g(C_t)$$
does not depend on $t\in \Delta$ and the claim follows.
\end{proof}

\medskip

\subsection{Foliations of general type}
It remains to study  the invariance of plurigenera for foliations of
general type.
The main ingredients used in the proof are Theorem \ref{contrEGL}, Theorem \ref{zeroloc},  Kawamata-Viehweg vanishing theorem and the study of the Zariski decomposition of $\lceil mP_t\rceil$.
\begin{prop}\label{p_kod2}
Let $(X_t,\mathcal{F}_t)_{t\in\Delta}$ be a family of foliations such that $\kappa(\mathcal{F}_t)=2$ for any $t\in\Delta$.
Then for any sufficiently large positive integer $m$, the dimension $h^0(X_t,mK_{\mathcal{F}_t})$ does not depend on $t\in \Delta$.
\end{prop}
\begin{proof} In the course of the proof, we denote by $\{D\}$ the fractional part of a $\mathbb Q$-divisor $D$.  
We divide the proof in 6 Steps: 

\smallskip
\noindent \textbf{Step 1.}
By Lemma \ref{Fecc}, after possibly shrinking $\Delta$, we may assume that $\mathcal F_0$ is relatively minimal,  where $0\in \Delta$ corresponds to the central fibre $X_0$.
By Proposition \ref{zardecfol}, we may write
$$K_{\mathcal{F}}=P+N$$ 
such that for any $t\in \Delta$, if we denote $P_t=P|_{X_t}$ and $N_t=N|_{X_t}$ then 
$$K_{\mathcal{F}_t}=P_t+N_t$$
is a Zariski decomposition of $K_{\mathcal F_t}$. After possibly shrinking $\Delta$ further, by Proposition \ref{zardecfol}
we may assume that each irreducible component of $N$ meet the surfaces $X_t$ transversally in a rational curve.
In particular, we have
$$H^0(X_t,mK_{\mathcal{F}_t})\cong H^0(X_t,\lceil mP_t\rceil),$$
for any $t\in \Delta$. It is enough to show  that 
$h^0(X_t,\lceil mP_t\rceil)\ge h^0(X_0,\lceil mP_0\rceil)$ for any $t\in \Delta$.
We denote by $i$ the Cartier index of $P_t$. By Lemma \ref{lcm}, $i$ does not depend on $t\in \Delta$.

\smallskip 
\noindent \textbf{Step 2.} Let $E_1,\dots, E_k$ be the irreducible components of $N$ and let 
$$E_i^t={E_i}|_{X_t}.$$ 
By Theorem \ref{structneg}, $E_i^t$ is a smooth rational curve. In particular, the self-intersection $(E_{i}^t)^2$ is a negative number which does not depend on $t\in \Delta$. 

Thus, if  $$\nu_t\colon X_t\rightarrow Y_t$$
is the contraction of ${\rm Supp}N_t$, 
then there exists an effective $\mathbb{Q}$-divisor $\Theta$ on $\mathcal{X}$, whose coefficients are contained in $(0,1)\cap\mathbb{Q}$ and whose support is contained in the support of $N$, such that
 $$K_{X_t}=\nu_t^{\ast}K_{Y_t}-\Theta\vert_{X_t},$$
for any $t\in \Delta$. We denote 
$$\Theta_t=\Theta\vert_{X_t}$$ for any $t\in \Delta$. 

\smallskip

\noindent \textbf{Step 3.}
For any $t\in \Delta$ and for any positive integer $m$, let $$\lceil mP_t\rceil=P^{(t)}_m+N^{(t)}_m$$
be the Zariski decomposition of $\lceil mP_t\rceil$.
Note that 
$$\lceil mP_t\rceil= mP_t+ \{-mP_t\}$$
and therefore $N_m^{(t)}\le  \{-mP_t\}$.
In particular,  the support of $N^{(t)}_m$  is contained in  the support of $N_t$ and
$$\lfloor N_m^{(t)}\rfloor=0$$
for any $t\in \Delta$.
Note that the divisor $N^{(t)}_m$ is uniquely determined by the intersection
of $\lceil mP_t\rceil$ with the components of ${\rm Supp}N_t$.
Thus, 
$$N^{(t)}_m=N^{(t)}_{m'}$$
for any  positive integers $m,m'$ which are equal modulo $i$. 
In addition,  there exist $\mathbb{Q}$-divisors $N_m$, $P_m$
on $\mathcal{X}$ such that $$\lceil mP\rceil=P_m+N_m$$
and $$N_m\vert_{X_t}=N^{(t)}_m,$$
for any $t\in \Delta$.

\smallskip
\noindent\textbf{Step 4.}
Let $\Gamma_1,\dots,\Gamma_p$ be all the elliptic Gorenstein leaves contained in the central fibre $X_0$ (cf. Definition \ref{d_egl}). By Theorem \ref{zeroloc}, there exist $\mathcal{C}_1,\dots,\mathcal C_q$ disconnected chains of rational curves in $X_0$ 
such that $$\{C\subseteq X_0\mid P\cdot C=0\}=\bigcup_{i=1}^p \Gamma_i\cup\bigcup_{i=1}^q \mathcal{C}_i\cup{\rm Supp}N_0.$$
Let $Z=\sum_{i=1}^p \Gamma_i+\sum_{i=1}^q \mathcal{C}_i$.

By Remark \ref{noint},   the curves $\Gamma_1,\dots,\Gamma_p$ do  not intersect the support of $N_0$ and  in particular $$\mathcal{O}_{\cup\Gamma_i}(\lceil mP_0 \rceil)=\mathcal{O}_{\cup\Gamma_i}(mK_{\mathcal F_0})$$
Thus, Theorem \ref{contrEGL} implies that for any positive integer $m$ the sheaf  $\mathcal{O}_{\cup\Gamma_i}(\lceil mP_0 \rceil)$ has degree zero and is not torsion. For each $i=1,\dots,p$, the curve $\Gamma_i$ is Cohen-Macaulay with trivial dualizing sheaf and  Serre duality implies that 
$$h^1(\Gamma_i, \lceil mP_0 \rceil)=h^0(\Gamma_i,-\lceil mP_0 \rceil)=0.$$
Thus 
\begin{equation}\label{e_gammai}
h^0(\cup\Gamma_i,\lceil mP_0 \rceil)=h^1(\cup\Gamma_i,\lceil mP_0 \rceil)=0.
\end{equation}
Let $C\subseteq \cup\mathcal{C}_i$ be an irreducible component.
By Theorem \ref{zeroloc}, there are two possibilities:
\begin{description}
\item[type A] $C\cap {\rm Supp}N_0=\emptyset$; or
\item[type B] $N_0\vert_C=\frac 12 p_1+\frac 12 p_2$, with $p_1,p_2\in C$ distinct points.
\end{description}
Moreover, if $C$ is of type B, then
$p_i$ belong to a connected component $E^C_{i}$ of ${\rm Supp}N_0$
which is a smooth rational curve of self-intersection -2, for $i=1,2$.
In particular, the coefficient of $\Theta_0$ along $E^C_i$ is zero for  $i=1,2$.

For any positive integer $m$, we have 
$$\deg \lceil mP_0\rceil\vert_C=mK_{\mathcal F_0}\cdot C- \deg \lfloor mN_0\vert_C\rfloor=
\begin{cases}
0& \quad{\rm if}\;C\;{\rm is}\;{\rm of}\;{\rm type}\;{\rm A}\\
0& \quad {\rm if}\;C\;{\rm is}\;{\rm of}\;{\rm type}\;{\rm B}\;{\rm and}\;m\;{\rm is}\;{\rm even}\\
1& \quad {\rm if}\;C\;{\rm is}\;{\rm of}\;{\rm type}\;{\rm B}\;{\rm and}\;m\;{\rm is}\;{\rm odd}.\\
\end{cases}
$$

Let $h_i$ be positive integers so that $$\lfloor N^{(0)}_m +\Theta_0 +\frac 12 \sum( E^C_{1}+E^C_{2}) \rfloor=\sum h_i E^0_i$$
where the sum runs over all the curves $C\subseteq \cup \mathcal C_i$ of type B. We first prove the following: 

\medskip

\noindent \underline{Claim:} We have $h^1(\cup\mathcal{C}_i, \lceil mP_0 \rceil-\sum h_i E^0_i)=0$.

\medskip

Let $C_0$ be an irreducible component of $\cup \mathcal C_i$. If $C_0$ is of type A, then $$\deg (\lceil mP_0 \rceil-\sum h_i E^0_i)\vert_{C_0}=0, $$
and in particular  $h^1(C_0, \lceil mP_0 \rceil-\sum h_i E^0_i)=0$.

If $C_0$ is of type B and $m$ is even, then the coefficient of $mP_0$ along $E^{C_0}_{i}$ is integral and in particular the coefficient of $N_m^{(0)}$ along $E_i^{C_0}$ is zero.
Since also the coefficient of $\Theta_0$ along $E^{C_0}_i$ is zero, 
we have that 
 $${\rm coeff}_{E_{i}^{C_0}}(\lfloor N^{(0)}_m +\Theta_0 +\frac 1 2\sum ( E^C_{1}+E_{2}^C) \rfloor)=0\qquad\text{for }i=1,2$$
and, in particular, $$\deg (\lceil mP_0 \rceil-\sum h_i E^0_i)\vert_{C_0}=0.$$
Thus, also in this case, we have $h^1(C_0, \lceil mP_0 \rceil-\sum h_i E^0_i)=0$.

Finally, if $C_0$ is of type $B$ and $m$ is odd, similarly as above we have
$$\deg(\sum h_i E^0_i)\vert_{C_0}\in\{0,1,2\},$$
and in particular $$\deg(\lceil mP_0 \rceil-\sum h_i E^0_i)\vert_{C_0}\in\{1,0,-1\}$$
and also in this case, we have $h^1(C_0, \lceil mP_0 \rceil-\sum h_i E^0_i)=0$.

Thus, if $\mathcal C_i = C_1\cup \dots\cup C_q$ is a connected component of $\cup \mathcal C_i$, then it is a chain of rational curves which admits at most one  tail  of type $B$, and all the other curves are of type $A$. Thus we have a short exact sequence 
$$0\to \mathcal O_{\mathcal C_i}(\lceil mP_0 \rceil-\sum h_i E^0_i)
\to \bigoplus_{j=1}^q\mathcal O_{C_j}(\lceil mP_0 \rceil-\sum h_i E^0_i) 
\to \bigoplus_{j=1}^{q-1}\mathcal O_{C_j\cap C_{j+1}} (\lceil mP_0 \rceil-\sum h_i E^0_i)\to 0$$
such that the induced map
$$\bigoplus_{j=1}^qH^0(C_j,\lceil mP_0 \rceil-\sum h_i E^0_i) 
\to \bigoplus_{j=1}^{q-1} H^0(C_j\cap C_{j+1}, \lceil mP_0 \rceil-\sum h_i E^0_i)$$ is surjective. 
Since $h^1(C_j, \lceil mP_0 \rceil-\sum h_i E^0_i)=0$ for any irreducible component $C_j$ of $\mathcal C$, 
the claim follows. 

\medskip 

Consider now the short exact sequence
$$0\rightarrow\mathcal{O}_{X_0}\left(\lceil mP_0 \rceil-Z-\sum h_i E^0_i\right)
\rightarrow\mathcal{O}_{X_0}\left(\lceil mP_0 \rceil-\sum h_i E^0_i\right)
\rightarrow\mathcal{O}_{\cup\Gamma_i\cup \mathcal{C}_i}\left(\lceil mP_0 \rceil-\sum h_i E^0_i\right)\rightarrow0.$$

Thus, \eqref{e_gammai} and the claim above yield the exact sequence 
$$
\begin{array}{rll}
0\rightarrow & H^0\left(X_0,\lceil mP_0 \rceil-Z-\sum h_i E^0_i\right)\rightarrow&
H^0\left(X_0,\lceil mP_0 \rceil-\sum h_i E^0_i\right)\rightarrow H^0\left(\cup\mathcal{C}_i,\lceil mP_0 \rceil-\sum h_i E^0_i\right)\rightarrow\\
\rightarrow & H^1\left(X_0,\lceil mP_0 \rceil-Z-\sum h_i E^0_i\right)\rightarrow&
H^1\left(X_0,\lceil mP_0 \rceil-\sum h_i E^0_i\right)\rightarrow0
\end{array}
$$
and
$$H^2\left(X_0,\lceil mP_0 \rceil-Z-\sum h_i E^0_i\right)
\cong H^2\left(X_0,\lceil mP_0 \rceil-\sum h_i E^0_i\right).$$

\smallskip
\noindent \textbf{Step 5.}
Set $\bar{P}_0=\nu_{0\ast}P_0$. By the Negativity Lemma, it follows that $P_0=\nu_0^\ast\bar{P}_0$. 

Let $a$ be a positive integer. By Step 3,  we have that $${\rm Supp}N^{(0)}_a\subseteq {\rm Supp}N_0.$$ 

Set $\bar{Z}=\nu_{0\ast}Z$. Thus $$Z=\nu_0^{\ast}\bar{Z}-\frac{1}{2}\sum (E^C_1+E_{2}^C),$$
where, as above, the sum is over all the curve $C$ of type $B$.
Let  $C_0\subseteq X_0$ be an irreducible curve such that $P_0\cdot C_0=0$.  
If $C_0$ is contained in the support of $N_0$, then 
$$(K_{X_0}+\Theta_0)\cdot C_0=0\qquad \text{and}\qquad (Z+ \frac 1 2 \sum (E^C_1+E_{2}^C))\cdot C_0=0.$$
On the other hand, if $C_0$ is a smooth rational curve which is not contained in the support of $N_0$, then Theorem \ref{zeroloc} and Remark \ref{noint} imply that 
$$(\Theta_0+Z+ \frac 1 2\sum (E^C_1+E_{2}^C))\cdot C_0\le C_0^2 +2.$$
Finally, if $C_0$ is a single rational nodal curve, then 
$$(K_{X_0}+Z)\cdot C_0=0\qquad \text{and}\qquad \Theta_0+ \frac  1 2\sum (E^C_1+E_{2}^C))\cdot C_0=0.$$
Thus, Theorem \ref{zeroloc} implies that for each  $C_0\subseteq X_0$  irreducible curve such that $P_0\cdot C_0=0$, we have
$$(K_{Y_0}+\bar{Z})\cdot \nu_{0\ast}C_0=(K_{X_0}+\Theta_0+Z+ \frac  1 2\sum (E^C_1+E_{2}^C))\cdot C_0\leq 0.$$ 
It follows that there exists a sufficiently large positive integer $b$ such that
$$bi \bar{P}_0-(K_{Y_0}+\bar{Z})$$ is big and nef. 

We denote $m=a+bi$. 
In particular,  by Step 3 we have that 
$N^{(0)}_a=N^{(0)}_m$.
In addition, 
$$\lceil mP_0 \rceil=biP_0+\lceil aP_0 \rceil.$$
 We have
$$
\begin{array}{rl}
&\lceil mP_0 \rceil-\sum \Gamma_i-\sum\mathcal{C}_i-\sum h_i E^0_i\\
=&\lceil mP_0 \rceil-Z-\sum h_i E^0_i\\
=& K_{X_0}+\nu_0^{\ast}(bi \bar{P}_0-(K_{Y_0}+\bar{Z}))+
\lceil aP_0 \rceil-\sum h_i E^0_i+\Theta_0+\frac{1}{2}\sum (E_{1}^C+E_{2}^C)\\
=& K_{X_0}+\nu_0^{\ast}(bi \bar{P}_0-(K_{Y_0}+\bar{Z}))+
P_a^{(0)}+N_a^{(0)}-\sum h_i E^0_i+\Theta_0+\frac{1}{2}\sum (E_{1}^C+E_{2}^C).
\end{array}
$$
 Since
 $\sum h_i E^0_i   =\lfloor N^{(0)}_m +\Theta_0 +\frac 12 \sum( E^C_{1}+E^C_{2}) \rfloor$, the divisor 
 $$N_a^{(0)}-\sum h_i E^0_i+\Theta_0+\frac{1}{2}\sum (E_{1}^C+E_{2}^C)$$ is an effective divisor 
 with coefficients in the interval $(0,1)$ and whose support is contained in the support of $N_0$. 
 On the other hand,
  $ \nu_0^{\ast}(bi \bar{P}_0-(K_{Y_0}+\bar{Z}))+P^{(0)}_a $ is big and nef. 
 
 Thus, Kawamata-Viehweg vanishing theorem implies 
 $$H^j(X_0,\lceil mP_0 \rceil-Z-\sum h_i E^0_i)=0\qquad \text{for all } j>0.$$
Then, by Step 4, we have
$$H^j(X_0,\lceil mP_0 \rceil-\sum h_i E^0_i)=0\qquad \text{for all } j>0.$$
In particular, $$h^0(X_0,\lceil mP_0 \rceil-\sum h_i E^0_i)=
\chi(X_0,\lceil mP_0 \rceil-\sum h_i E^0_i).$$

Note that  $\chi(X_t,\lceil mP_t \rceil-\sum h_i E^t_i)$ does not depend on $t\in \Delta$. 

\medskip 
\noindent\textbf{Step 6.}
In Step 4 we defined
$\sum h_i E^0_i =\lceil N^{(0)}_m +\Theta_0+\frac 1 2\sum (E_{1}^C+E_{2}^C) \rceil$.
Note that, since the coefficient of $\Theta_0$ are contained in the interval $(0,1)$, and for each curve $C$ of type $B$, the curves $E_1^C$ and $E_2^C$ are not contained in the support of $\Theta_0$,  it follows that 
$$\lfloor N^{(0)}_m +\Theta_0+\frac 1 2\sum  (E_{1}^C+E_{2}^C) \rfloor\leq\lceil N_m^{(0)} \rceil.$$ 
Thus, we have that 
$$\sum h_i E^t_i \le \lceil N_m^{(t)}\rceil
$$
for any $t\in \Delta$. 
By the properties of the Zariski decomposition,
$$H^0(X_t,\lceil mP_t \rceil)\cong H^0(X_t,\lfloor\lceil mP_t \rceil-N_m^{(t)}\rfloor)
=H^0(X_t,\lceil mP_t \rceil-\lceil N^{(t)}_m \rceil).$$
More precisely, for any effective divisor $E\leq \lceil N^{(t)}_m \rceil$,
$$H^0(X_t,\lceil mP_t \rceil)\cong H^0(X_t,\lceil mP_t \rceil-E).$$
In particular
$$h^0(X_t,\lceil mP_t\rceil)=h^0(X_t,\lceil mP_t \rceil-\sum h_i E^t_i)$$
for all $t\in \Delta$. 
Note that since $P_t$ is big and nef for all $t\in \Delta$, if $m$ is sufficiently large then Serre duality implies that
$$h^2 (X_t,\lceil mP_t \rceil-\sum h_i E^t_i)=0$$ 
for all $t\in \Delta$. 
Since $$\lceil mP_t \rceil-\sum h_i E^t_i=(\lceil mP \rceil-\sum h_i E_i)|_{X_t}$$
for all $t\in \Delta$, Step 5 implies that 
$$\begin{aligned}
h^0(X_t,\lceil mP_t \rceil)&= h^0(X_t,\lceil mP_t \rceil-\sum h_i E^t_i)\\
&\ge \chi(X_t,\lceil mP_t \rceil-\sum h_i E^t_i) \\
&= \chi(X_0,\lceil mP_0 \rceil-\sum h_i E^0_i) \\
&=h^0(X_0,\lceil mP_0 \rceil-\sum h_i E^0_i)=h^0(X_0,\lceil mP_0 \rceil)
\end{aligned}$$
for all $t\in \Delta$, concluding the proof of the Proposition.
\end{proof}
\medskip

We are now ready to prove our main Theorem:
\begin{proof}[Proof of Theorem \ref{t_main}]
Proposition \ref{p_kod0} implies (1). 
Proposition \ref{p_kod1} together with Propositions  \ref{p_kod1a}, \ref{p_kod1b} and \ref{p_kod1d} imply $(2)$  and $(3)$.  
Finally, Proposition \ref{p_kod2} implies $(4)$. 
\end{proof}

\medskip

\subsection{Some Examples}
We now show three examples on which invariance of plurigenera does not hold. 

We begin by providing an example of a family of foliations which does not satisfy 
hypothesis (2) of Definition \ref{hyp} and for which 
the Kodaira dimension is not constant: 
\begin{ex}\label{e}
For $j=1,2$, let $f_{j, t_0}\in \mathbb C[x_0,x_1,x_2]$ be a homogeneous polynomial of degree four which is the product of
a linear factor $l$ and a factor $c_j$ of degree 3 such that $\{c_j=0\}$ is a smooth cubic
$$f_{j, t_0}=l\cdot c_j.$$
For $j=1,2$, let $f_{j, t_1}\in \mathbb C[x_0,x_1,x_2]$ be a homogeneous polynomial of degree four such that $\{f_{j, t_1}=0\}$
is a smooth quartic.
We may assume that $f_{j, t_0}$ and $f_{j, t_1}$ are such that
\begin{itemize}
\item for any $t\in \Delta\setminus \{0\}$ and $j=1,2$, if we define  $f_{j,t} = t f_{j, t_0}+(1-t)f_{j, t_1}$ then the curve
$\{f_{j,t}=0\}$ is a smooth quadric;\item for any $t\neq 0$ the curves $\{f_{1,t}=0\}$ and $\{f_{2,t}=0\}$ meet transversally;
\item the curves $\{c_1=0\}$ and $\{c_2=0\}$ meet transversally in 9 points $p_1,\dots,p_9$; and
\item for any $t\in \Delta\setminus\{ 0\}$, the curves of the pencil 
$$\mathcal{P}_t=\{uf_{1,t}+vf_{2,t}=0\left\vert\; [u:v]\in\mathbb{P}^1\right.\}$$
are all irreducible and reduced.
\end{itemize}
The base points of the pencils $\mathcal{P}_t$ form 16 curves $\mathcal{B}_1,\dots,\mathcal{B}_{16}$ in $\mathbb P^2\times \Delta$
meeting the fibres transversally.
After possibly reordering the points $p_1,\dots,p_9$, we may assume that there exists $k$ such that the curves  $\mathcal{B}_1,\dots,\mathcal{B}_{16}$ pass through $p_{k+1},\dots,p_9$
and we may pick $\mathcal{B}_{17},\dots,\mathcal{B}_{16+k}$
 smooth curves in $\mathbb P^2\times \Delta$ meeting $\mathbb{P}^2\times \{0\}$ transversally in $p_1,\dots,p_k$.

Let  $$\mathcal{X}\xrightarrow{\varepsilon} \mathbb{P}^2\times\Delta\rightarrow\Delta$$
be the blow-up of $\mathbb{P}^2\times\Delta$ along the curves $\mathcal{B}_1,\dots,\mathcal B_{16+k}$. Let 
$\pi\colon \mathcal X\to \Delta$ be the induced morphism with $X_t=\pi^{-1}(t)$ and let 
$\varepsilon_t\colon X_t\to \mathbb P^2$ be the induced morphism with exceptional divisors $E^t_1,\dots,E_{16+k}^t$.

For any $t\in\Delta\backslash\{0\}$ we have a fibration $f_t\colon X_t\rightarrow\mathbb{P}^1$
whose general fibres are the strict transforms of the elements of $\mathcal{P}_t$.
On $X_0$ we have a fibration $f_0\colon X_0\rightarrow\mathbb{P}^1$
whose fibres are the curves in the pencil generated by $c_1$ and $c_2$.
Let $H$ be a hyperplane section in $\mathbb{P}^2$ and let $p\in \mathbb P^1$ be a general point.
For any $t\neq 0$, the fibre $f_t^{\ast}p$ is the strict transform of a curve in the pencil,
thus $$f_t^{\ast}p\sim \varepsilon_t^{\ast}(4H)-\sum_{i=1}^{16}E^t_i.$$ 

Let $\mathcal F_t$ be the foliation on $X_t$ induced by the fibration $f_t$. 
Then for any $t\neq 0$
$$K_{\mathcal{F}_t}=K_{X_t/\mathbb{P}^1}=
\varepsilon_t^{\ast}K_{\mathbb{P}^2}+
\sum_{i=1}^{16+k} E^t_i+
f_t^{\ast}(2p)\sim \sum_{i=17}^{16+k}E^t_i+\varepsilon^{\ast}_t H+f_t^{\ast}p.$$
In particular, $K_{\mathcal F_t}$  is  big.
On the other hand, if $t=0$  then
$$f_0^{\ast}p\sim \varepsilon_0^{\ast}(3H)-\sum_{i=1}^{9}E^0_i.$$
and
$$K_{\mathcal{F}_0}=K_{X_0/\mathbb{P}^1}=
\varepsilon_0^{\ast}K_{\mathbb{P}^2}+
\sum_{i=1}^{16+k} E^0_i+ f_0^{\ast}(2p)\sim \sum_{i=10}^{16+k} E^0_i+ f_0^{\ast}p.$$
Then for $t=0$ the canonical divisor $K_{\mathcal{F}_0}$ is not big.

Note that the foliation we obtain is induced by a fibration on each fibre,
but there does not exist a fibration  on $\mathcal X$ which induces the foliation.
\end{ex}

\medskip

We now show an example of a family of foliations $(X_t,\mathcal F_t)_{t\in\Delta}$ of Kodaira dimension one, which are induced by nonisotrivial elliptic fibrations and such that $h^0(X_t,\mathcal O_{X_t}(K_{\mathcal{F}_t}))$
is not constant, for $t\in \Delta$.

\begin{ex}\label{e_elliptic}
Let $C$ be a curve of genus at least 2.
Let $p_{0,t}$, $p_{1,t}$ be two families of points on $C$, with $t\in \Delta$, 
such that $p_{0,0}=p_{1,0}$ and $p_{0,t}\neq p_{1,t}$ for any $ t\neq 0$.
The line bundle $$\mathcal{L}_t=\mathcal{O}_C(K_C+p_{0,t}- p_{1,t})$$
defines an elliptic fibration
$$\varphi_t\colon X_t\rightarrow C$$
such that its moduli part is $M_{C,t}=K_C+p_{0,t}- p_{1,t}$ and with discriminant $B_{C,t}=0$ for any $t\in \Delta$ (cf. \eqref{e_disc}).
Let $\mathcal{F}_t$ be the foliation associated to $\varphi_t$. By Remark \ref{r_redsnc}, the singularities of $\mathcal F_t$ are reduced and therefore they define a family of foliations  $(X_t,\mathcal{F}_t)_{t\in\Delta}$.
We have $$H^0(X_t,\mathcal O_{X_t}(K_{\mathcal{F}_t}))=H^0(C,\mathcal O_{C_t}(M_{C,t}))$$
and the latter vector space has dimension $g$ if $t=0$ and strictly less than $g$ if $t\neq 0$.
\end{ex}

\medskip

Note that elliptic Gorenstein leaves (cf. Definition \ref{d_egl}) never appear on foliations of general type induced by fibrations over a curve. Indeed, 
let $\mathcal F$ be a foliation of general type on a smooth surface $X$ induced by a fibration $f\colon X\to C$ over a curve $C$. 
Assume by contradiction that there exists an e.g.l. $\Gamma$ on $X$. Since $
\Gamma$ is $\mathcal{F}$-invariant,
$\Gamma$ is contained in a fibre $F$ of $f$.
Since the fibres of $f$ are connected, Remark \ref{noint} implies  that ${\rm Supp}(\Gamma)={\rm Supp}(F)$.
Since $K_{\mathcal F}\cdot \Gamma_i=0$ for any curve contained in the support of $\Gamma$, it follows  that $K_{\mathcal{F}}\cdot F=0$,
which is a contradiction because $K_{\mathcal{F}}$ is big and the fibres form a covering family.

Nevertheless, even for foliations of general type induced by fibrations,
invariance of plurigenera does not hold for small values of $m$, as our example below shows. 
\begin{ex}

Let $C$ be a smooth curve of genus $g\ge 2$. Let $\mathcal L$ be a line bundle on $\mathcal C=C\times \Delta$ such that, if we denote $L_t=c_1(\mathcal L|_{C\times \{t\}})$, then 
$$L_0=K_C\qquad\text{and}\qquad h^0(C,L_t)<h^0(C,K_C)=g\quad \text{ if $t\neq 0$}.$$
Let  $$\mathcal{E}=\mathcal{O}_{\mathcal C} \oplus \mathcal{O}_{\mathcal C} \oplus \mathcal L$$
and let $p\colon \mathcal Z=\mathbb P(\mathcal E)\to\mathcal  C$. Let $\eta \colon \mathcal Z \to \Delta$ be the induced morphism. For any $t\in \Delta$, we denote by $p_t\colon Z_t:=\eta^{-1}(t)\to C$  the restriction morphism. 
We have
$$K_{Z_t/C}=-3\xi_t + p_t^{\ast}L_t$$
where $\xi_t =c_1(\mathcal{O}_{\mathcal E_t}(1))$. Let $\xi=c_1(\mathcal O_{\mathcal E}(1))$.  The linear system $|4\xi|$ is base point free and in particular the general element 
$\mathcal X\in |4\xi|$ is smooth. Let $\pi\colon \mathcal X\to \Delta$ be the induced morphism and let $X_t=\pi^{-1}(t)$. Then, the 
 morphism $\mathcal X\to\mathcal C$ defines a regular foliation $\mathcal F$ on $\mathcal X$ such that if $\mathcal F_t$ is the restriction of $\mathcal F$ to $X_t$ then \eqref{folfibr} implies that
$$K_{\mathcal F_t}=K_{X_t/C}=(\xi_t+ p_t^{\ast}L_t)|_{X_t}$$
for all $t\in \Delta$. Note that $\mathcal F_t$ is a foliation of general type. 
We want to show that $h^0(X_t,K_{\mathcal F_t})$ is not constant. 

The dimension
$$h^0(Z_t,\xi_t+ p_t^{\ast}L_t)=2h^0(C,L_t)+h^0(C,2L_t)$$
is not constant by our choice of $\mathcal L$. 
Therefore. it is enough to show that 
$$h^0(X_t,K_{X_t/C})=h^0(Z_t,\xi_t+ p_t^{\ast}L_t).$$

Pick $t\in \Delta$. By the exact sequence obtained by restriction, we have

$$0\to H^0(Z_t,\xi_t+ p_t^{\ast}L_t)\to H^0(X_t,K_{X_t/C})\to H^1(Z_t,K_{Z_t/C}).$$
On the other hand, by Serre duality
$$h^1(Z_t,K_{Z_t/C})=h^2(Z_t,p^{\ast}_t K_C)$$
and the latter dimension is zero by the Leray spectral sequence. Thus, the claim follows. 

\end{ex}

\medskip

\bibliographystyle{amsalpha}
\bibliography{Library}

\end{document}